\documentclass[]{article}
\RequirePackage{epsfig}
\RequirePackage{amsbsy}
\RequirePackage{amsfonts}
\RequirePackage{amsmath}
\RequirePackage{amssymb}

\newtheorem{theorem}{Theorem}[section]
\newtheorem{lemma}{Lemma}[section]
\newtheorem{proposition}{Proposition}[section]
\newtheorem{corollary}{Corollary}[section]
\newtheorem{definition}{Definition}[section]
\newtheorem{example}{Example}[section]
\newtheorem{remark}{\bf Remark}[section]

\begin{document}
\title{{\itshape Test on the components of mixture densities}}

\author{Florent AUTIN\footnote{Address : C.M.I., 39 rue F. Joliot Curie, 13453 Marseille Cedex 13. Universit\'e Aix-Marseille 1. FRANCE. Email: autin@cmi.univ-mrs.fr}\hspace{2pt}
(Universit\'e Aix-Marseille 1) \\
and Christophe POUET\footnote{Address : Ecole Centrale Marseille, 38 rue F. Joliot-Curie, 13451 Marseille Cedex 20. FRANCE. Email : cpouet@centrale-marseille.fr}
\hspace{2pt} (Ecole Centrale Marseille)}

\maketitle

\begin{abstract}
This paper deals with statistical tests on the components of mixture densities.  We propose to test whether the densities of two
 independent samples of independent random variables 
$Y_1, \dots, Y_n$ and $Z_1, \dots, Z_n$ result from the same mixture of $M$ components or not. We provide a test procedure which 
is proved to be asymptotically optimal according
to the minimax setting. We extensively discuss 
the connection between the mixing weights and the performance of the
 testing procedure and illustrate it with numerical examples. 
This link had never been clearly exposed up to now.
 \end{abstract}
 

\renewcommand{\thefootnote}{}
\footnotetext{\hspace*{-.51cm}AMS 2010 subject classification: Primary: 62C20, 62G10, 62G20; Secondary: 30H25, 42C40 \\ %
}

\section{Introduction}
\subsection{Mixture model with varying mixing weights}
Since more than 20 years, the mixture model has gained a lot of attention. This is due to its ease of interpretation by viewing each 
component as a distinct group in the data. This model has been
widely applied in several areas such as finance, economy, biology,
astronomy, survey methods,...

\noindent  Most of the theoretical results in the literature deal with
the estimation of the components or of the mixing weights. 
There are two types of mixture models : the most popular one has fixed mixing weights and the other one has varying mixing weights.\\
On the one hand, many statisticians have been interested
in estimating the mixing weights. For example, Hall \cite{Hall}, Titterington \cite{Tit} and Hall and Titterington \cite{HT-1} have considered nonparametric estimation
of the mixing weights. Two other examples
about the mixing weights
are the estimation of a functional of the weights by van de Geer \cite{vdG} and 
the computation of confidence intervals by Qin \cite{Q}.
On the other hand, one can be interested in estimating the components of
the mixture. This can be easily 
done with varying mixture weights by applying several well-known 
methods such as histograms in Lodakto and Maiboroda \cite{LM-1}, empirical distribution
in Maiboroda \cite{Mai-3} or wavelet thresholding methods in Pokhyl'ko \cite{Pok-1}. Finally, the mixing weights and the mixture components can also 
be estimated both and at the same time, this result holds in a particular
setting for k-variate data introduced by Hall and Zhou \cite{HZ}.\\

\noindent More recently, the mixture model has also been studied in the testing problem framework.
The usual addressed question is whether
 the observations come from a non-trivial
 mixture model or from a trivial one (i.e. with
only one component). This has been done for example by 
Garel \cite{Gar-1} and \cite{Gar-2} and Delmas \cite{Del} in the case of fixed
mixing weights and by Maiboroda \cite{Mai-2} in the case of
 varying mixing weights.
Their homogeneity tests which rely respectively on the likelihood ratio test and on a Kolmogorov-Smirnov type test are proved to be consistent.
Here we propose to study a testing problem with two samples
in a mixture model with varying mixing weights. \\
\noindent Although the varying mixing weights model does not seem natural at first sight, on can think of several 
situations where it can be useful. Let us give three examples that will help the reader
to recognize its usefulness.
\medskip

\noindent \underline{Social science}

\noindent This first example is the closest to the varying mixing weights model that is studied here.
Let us consider an organization divided into several departments such as an enterprise. Aggregated informations are only
known at the department level, e.g. proportion of men and women, proportion of graduates and undergraduates,
proportion of married and unmarried people,etc... The researcher is interested in a variable for these subgroups such as salary. For each person, the researcher has only recorded salary and department. The information of
interest which allows to divide the sample into subgroups is unavailable
at the individual level.
 This can happen if the 
researcher has forgotten to record this information when collecting the data; this frequently happens when a new question
arises during the study of the data. Another reason can be that the law forbids to record
such information at the individual level; for example this is the case of origins or races in many countries.
There is a wealth of works on partially missing data (see
McKnight et al. \cite{Knight} for example) but the case of entirely
missing data has never been really considered. 
From our point of view, a varying mixing
 weights model is a way to cope with
 this lack of information at the individual level and
to allow the researcher to reconstruct information for each subgroup. 
Although we are aware of methodological problems,
we want to emphasize that in this case the varying mixing weights are exactly known to the researcher; indeed, aggregated
information often exists and is much easier to collect than
individual information.
\newpage
\noindent \underline{Image analysis}

\noindent Let us assume a simple picture taken at a party and consisting of people and background. 
There is usually no way to distinguish at the pixel level whether it comes 
from people or background.
Nevertheless one can think of some kind of aggregated information to roughly divide the image
into several areas.
In the center of the image, there are
usually mainly people and only little background. In the area surrounding the center, there is mainly
background although few people can be scattered here and there. Therefore the image can be divided
into two areas. This written description of the image can be translated into a mathematical
description namely the varying mixing weights model. In this model, the statistician will
be able to extract distinctive features of the picture concerning people or 
background.
We are aware that methodological problems can appear in this setup. 
For example, spatial structure are not taken into account.
One can consider that the
weights are only roughly known which can be a problem. 
Nevertheless for some types of images, such as satellite
images, one can assume that the weights are accurately known. 
Indeed, as the area under scrutiny is exactly known from a geographical
point of view, one can use
aggregated information about surfaces such as proportion of forest, land, city, water,etc...
\medskip

\noindent \underline{Finance}

\noindent Mixture densities have been proved to be useful in volatility modeling
(see Bernhard and Leblang \cite{BL}, Avellaneda \cite{Ave} for example). If one consider
the volatility clustering effect (see Cont \cite{Cont} for example), one can roughly
divide time into periods where the proportions of high and low volatility 
are estimated. Indeed during each period it might be hard
to exactly label observations corresponding to low or high volatility.
Therefore the varying mixing weights model can
be considered and help to extract useful features of the 
mixture components. This case with estimated proportions is not solved 
here. Although it is beyond the scope of this paper, we briefly discuss it 
in Section \ref{section4}.
\medskip

\noindent Let us now come back to our testing problem with two samples
in a mixture model with varying mixing weights: let $Y_1, \dots, Y_n$ and $Z_1, \dots, Z_n$ be 
two independent $n$-samples of independent random variables. We propose to study in this paper
whether these two samples of random variables come from the same mixture of $M$ unknown densities $p_{u}$ ($1 \leq u \leq M$) or not. 
We assume that the mixing weights associated with each observation are available to the statistician.
In Butucea and Tribouley \cite{BT} some procedures are proposed to test if two $n$-samples of i.i.d. variables have common probability density.
Their setting is equivalent to the case $M=1$ in our mixture problem. 
Here the problem appears more complex since the two samples are not based on random variables with the same marginal densities. 
Our results show that there is no loss in the minimax rate compared
to the simpler case studied by Butucea and Tribouley \cite{BT}.
In Section \ref{section2} we provide an asymptotically minimax test which is based
on wavelet methods and we prove the dependence between the  mixing weights and 
the constants appearing in the definition of the minimax rate of testing. Until now this phenomenon has never been studied and
is extensively discussed in this paper.  In addition to our theoretical result some numerical experiments are given  in Section 
\ref{section3} in order to illustrate  the strong connection between the mixing weights and the performance of the test. As expected, our test performs very well for various mixture models. Sections  \ref{section4} and  \ref{section5} are respectively devoted to possible extensions of work  and to proofs of main results.\\

\noindent Here we introduce the wavelet 
framework that will be used.

\subsection{Wavelet framework}
\noindent We first recall  that wavelets have been often applied in different mathematical fields such as in approximation theory, in signal analysis and in statistics for instance. 
In particular,
many recent statistical works on estimation (see among others Autin \cite{Aut}, Donoho et al \cite{DJKP}, Cohen et al \cite{CDKP} ) and on hypothesis testing (see Spokoiny \cite{Spok}) use the wavelet setting to provide efficient estimators and tests.
There are many explanations for the huge interest of the wavelet setting. One of them is that wavelets bases are localized both in frequency and in time, contrary to the classical
Fourier basis which is only localized in frequency. As a consequence, the wavelet setting appears to be well adapted to describe local characteristics of a signal to be reconstructed.\\

\noindent Let $\phi$ and$\psi$
be  two compactly supported functions of $L_2(\mathbb{R})$ and denote
for all $j$ in $\mathbb{N}
$ and all $k$ in $\mathbb{Z} \mbox{ and all } x$ in $\mathbb{R}$,
$\phi_{jk}(x)=2^{^{j/2}}\phi(2^{^j}x-k)$ and  $\psi_{jk}(x)=2^{^{j/2}}\psi(2^{^j}x-k)$. \\

\noindent
Suppose that for any $j$ in $\mathbb{N}$:
\begin{itemize}
\item  $\{\phi_{jk}, \psi_{j'k}; \: j'\geq j; k \in \mathbb{Z}\}$ constitutes an
orthonormal basis of $L_{2}(\mathbb{R})$,
\item  $support(\phi)\cup support(\psi)\subset [-L,L[$ for some $L>0$.
\end{itemize}

\noindent Some most popular examples of such
bases, called compactly supported orthonormal wavelet bases, are given in Daubechies \cite{Daub}.
The function $\phi$ is called the scaling function and $\psi$ the associated wavelet.\\

\noindent
Any function $h$ in $L_2(\mathbb{R})$  can be represented as:
$$h(t)=\sum_{k\in
\mathbb{Z}}\alpha_{jk} \phi_{jk}(t) + \sum_{j'\geq j}\sum_{k\in
\mathbb{Z}}\beta_{jk}\psi_{j'k}(t)$$ \noindent where
$\forall j \in \mathbb{N}, \forall j' \geq j, \forall k \in \mathbb{Z}$:
\begin{itemize}
\item
$\alpha_{jk}=\begin{displaystyle}\int\end{displaystyle}_{I_{jk}}h(t)\phi_{jk}(t)dt$ \quad and \quad 
$\beta_{j'k}=\begin{displaystyle}\int\end{displaystyle}_{I_{j'k}}h(t)\psi_{j'k}(t)dt$,
\item $I_{jk}=\left\{x \in \mathbb{R}; -L\leq 2^{j}x-k < L\right \}=\left[\frac{k-L}{2^j},\frac{k+L}{2^j}\right[.$\\
\end{itemize}

\noindent Let us now describe the testing problem we focus on. \\

\subsection{Mathematical description of the testing problem}
\noindent Let $Y_1, \dots, Y_n$ be a 
sample of  independent random variables with unknown marginal densities
$$f_i(.)=\sum_{u=1}^M \omega_u(i) p_u(.), \quad 1\leq i \leq n,$$
\noindent and let $Z_1, \dots, Z_n$ be another sample of  independent random variables with unknown marginal densities
$$g_i(.)=\sum_{u=1}^M \sigma_u(i) q_u(.), \quad 1\leq i \leq n.$$
We also assume that the two samples are independent.

\noindent Here and in what follows, we suppose that the mixing weights $(\omega_u(i), \: 1 \leq u \leq M$, $1 \leq i \leq n)$ and $(\sigma_u(i), \: 1 \leq u \leq M, 1 \leq i \leq n)$ are known to the statistician and satisfy
\begin{itemize}
\item $\forall (u,i) \in \{1, \dots, M \}\times \{1, \dots, n\}, \: \min(\omega_u(i),\sigma_u(i))\geq 0$,
\item $\forall i \in  \{1, \dots, n\}, \: \begin{displaystyle}\sum_{u=1}^M \end{displaystyle}\omega_u(i)=\begin{displaystyle}\sum_{u=1}^M \end{displaystyle}\sigma_u(i)=1$,
\end{itemize} and are known by the statistician whereas the densities $ p_u$ and $q_u$ $(1 \leq u \leq M)$ are unknown.\\

\noindent Let us denote
$ \overrightarrow{p}= \left( p_1, \dots, p_M \right)$ and $
\overrightarrow{q}= \left( q_1, \dots, q_M \right).$ \\
\noindent 
We study in this paper a nonparametric procedure to test whether the samples result from the same mixture of densities.
Let $\mathcal{D}$ denote the set of all probability densities  with 
respect to the Lebesgue measure on $\mathbb{R}$. For any real number $R>0$, 
we define
$$ \Theta_0\left( R \right) = \left\{ 
\left( \overrightarrow{p}, \overrightarrow{q} \right) :
\forall u \in \{1, \dots, M\}, \quad p_u = q_u \in \mathcal{S}(R) \right\}$$
where $\mathcal{S}(R)=\mathcal{D}\cap\mathbb{L}_\infty(R)\cap
\mathbb{L}_2(R)$.
\\ We consider the following null hypothesis
$$\mathcal{H}_0: \quad \left( \overrightarrow{p}, \overrightarrow{q} \right) \in \Theta_0\left( R \right).$$
For a given $C>0$, we define
\begin{eqnarray*}
\Theta_1\left( R,C,n,s\right) & =  & \Big\{ 
\left( \overrightarrow{p}, \overrightarrow{q} \right) :
\forall u \in \{1, \dots, M\}, p_u-q_u \in \mathcal{B}^s_{2,\infty}(R), \\ 
& & \quad \exists u \in \{1, \dots, M\}, \left(p_u,q_u\right) \in \Lambda_n(R,C) \Big\},
\end{eqnarray*}
\noindent where  $\Lambda_n(R,C)=\left\{(p,q) \in (\mathcal{D}\cap\mathbb{L}_\infty(R))^2, \| p- q\|_2 \geq C r_n\right\},$ for a sequence $r_n$ tending to $0$ when $n$ goes to infinity and $\mathcal{B}^s_{2,\infty}(R)$ is the $R$-ball of a functional space defined below. We consider the following alternative
\begin{eqnarray*}
\mathcal{H}_1 &:&  \quad \left( \overrightarrow{p}, \overrightarrow{q} \right) \in \Theta_1\left( R,C,n,s\right).
\end{eqnarray*}

\noindent As usual in the nonparametric setting, we focus on
a large class of functions having some regularity so as to derive optimal properties. For the chosen wavelet basis, the space $\mathcal{B}^s_{2,\infty}(R)$ represents the  
$R$-ball of the so-called Besov body which is composed of all the functions $h \in L_2(\mathbb{R})$ for which the sequence of wavelet coefficients $(\alpha_{jk}, \ \beta_{j'k}, j \in \mathbb{N}, j' \geq j, k \in \mathbb{Z})$
 satisfies:$$ \sup_{j \in \mathbb{N}}2^{2js}\sum_{j' \geq j}\sum_{k \in \mathbb{Z}} \beta_{j'k}^2 \leq R.$$

\noindent {\underline{The minimax setting}}\\
\noindent In this paragraph we recall the minimax approach which is often used to evaluate the performances of testing procedures.
 Given the sum of the probability errors,
say $\gamma \in \left[0,\ 1 \right]$, we study the optimal separation rate $r_n$ between the null hypothesis and the alternative. This rate $r_n$ is the best possible rate separating 
at least one of the $M$ couples of density components $p_u$ and $q_u$.
It is  usually called {\it{the minimax rate}}. Let us recall the classical definition for the separation rate.\\

\begin{definition}
Let $0<\gamma<1$. We say that $r_n$ is the minimax rate separating $\mathcal{H}_0$ and $\mathcal{H}_1$ of our testing problem at level $\gamma$ if the two following statements are satisfied:

\begin{enumerate}
\item there exist a sequence of test procedures $\Delta^*_n$ and  a constant $C_\gamma$ such that 
{\small{\begin{eqnarray}\label{uptheo0}
\!\!\!\!\!\!\!\!  
\limsup_{n \to \infty}  \left( \sup_{(\overrightarrow{p},\overrightarrow{q}) \in \Theta_0(R)}\mathbb{P}_{\overrightarrow{p},\overrightarrow{q}}(\Delta^*_n=1) + \sup_{(\overrightarrow{p},\overrightarrow{q})  \in \Theta_1(R,C,n,s)} \mathbb{P}_{\overrightarrow{p},\overrightarrow{q}}(\Delta^*_n=0)\right) \leq \gamma
\end{eqnarray}}}
for all $C>C_\gamma$;
\item  there exists a constant $c_\gamma$ such that 
{\small{\begin{eqnarray}\label{lowtheo0}
\!\!\!\!\!\!\!\! 
\liminf_{n to \infty}  
\inf_{\Delta}
\left( \sup_{(\overrightarrow{p},\overrightarrow{q}) \in \Theta_0(R)}\mathbb{P}_{\overrightarrow{p},\overrightarrow{q}}(\Delta=1) + \sup_{(\overrightarrow{p},\overrightarrow{q})  \in \Theta_1(R,C,n,s)} \mathbb{P}_{\overrightarrow{p},\overrightarrow{q}}(\Delta=0)\right) > \gamma
\end{eqnarray}}}
for all $C<c_\gamma$, where the infimum is taken over all test procedures $\Delta$.\end{enumerate}
\end{definition}

\noindent {\underline{Hypothesis on the model}}\\
\noindent In our study we suppose that the mixing weights $(\omega_u(i), \: 1 \leq u \leq M, 1 \leq i \leq n)$ and $(\sigma_u(i), \:  1 \leq u \leq M, 1 \leq i \leq n)$ satisfy an added hypothesis. Let us denote by 
$\Omega=(\Omega)_{u,i}$ the matrix  with coefficients $\Omega_{u,i}=\omega_{u}(i)$ and $\Sigma=(\Sigma)_{u,i}$ the matrix with coefficients $\Sigma_{u,i}=\sigma_{u}(i).$ 
\begin{itemize}
\item[HYP-1] The smallest eigenvalues of the $(M \times M)$-matrices $\Gamma_n=\Omega \Omega^*$ and $\Gamma'_n=\Sigma \Sigma^*$ are both larger than or equal to $Kn$, with $0<K<1.$
\end{itemize} 

\noindent 
We recall the following proposition due to Maiboroda \cite{Mai-3}.
\begin{proposition}\label{maibo}
Suppose that the previous conditions are satisfied by the mixing weights $(\omega_u(i), 1 \leq u \leq M, 1 \leq i \leq n)$ and $(\sigma_u(i), 1 \leq u \leq M, 1 \leq i \leq n)$ associated with the model. Then, 
there exists a solution 
 of the two problems \\
$\left[\right.$ {\it{find}} $a_l=\{a_l(i), i=1, \dots, n\}$ such that $ <\omega_k,a_l>_n := \frac{1}{n}\begin{displaystyle}\sum_{i=1}^n \end{displaystyle}\omega_k(i)a_l(i)= \delta_{kl} \left. \right],$ \\
\noindent $\left[ \right.$ {\it{find}} $b_l=\{b_l(i), i=1, \dots, n\}$ such that $ <\sigma_k,b_l>_n := \frac{1}{n}\begin{displaystyle}\sum_{i=1}^n \end{displaystyle}\sigma_k(i)b_l(i)= \delta_{kl} \left. \right],$ 

\noindent where $\delta_{kl}$ is the Kronecker delta. According to HYP-$1$
 this solution satisfies
\begin{eqnarray}\label{hypflo1}
\sum_{l=1}^M <a_l,a_l>_n  := 
\frac{1}{n}\sum_{l=1}^M \sum_{i=1}^n a_l^2(i)  \leq \frac{M}{K},
\end{eqnarray}
\begin{eqnarray}\label{hypflo2}
\sum_{l=1}^M <b_l,b_l>_n  := 
\frac{1}{n}\sum_{l=1}^M \sum_{i=1}^n b_l^2(i)  \leq \frac{M}{K}.
\end{eqnarray}
\end{proposition}

\section{Nonparametric test procedure}\label{section2}

\noindent This paragraph deals with the case where the regularity $s$ of the Besov body that appears in $\mathcal{H}_1$ is known. From now on we denote by $a_{l}$ and $b_{l}$  the $n$-vectors which are the solutions of the two optimization problems appearing in Proposition \ref{maibo}. 
Let us describe the asymptotically minimax decision rule.

\subsection{Definition of the test procedure}

\noindent For each level parameter $j$, we define the test procedure $\Delta_j$ comparing the test statistic
$$T_j=\frac{1}{n^2}\sum_{l=1}^M\sum_k\sum_{i_1 \not= i_2}
\left[a_l(i_1)\phi_{jk}(Y_{i_1})-b_l(i_1)\phi_{jk}(Z_{i_1}) \right]
\left[a_l(i_2)\phi_{jk}(Y_{i_2})-b_l(i_2)\phi_{jk}(Z_{i_2}) \right]$$
 with a threshold value $t_n=t \ r_n^2$ where $t$ is a constant chosen later. We define
$$\Delta_j = \left\{
    \begin{array}{ll}
        1 & \mbox{ if }  T_j > t_n, \\
        0 & \mbox{ if }  T_j \leq t_n.
    \end{array}
\right.$$

\subsection{Properties of the test statistic}

\noindent In this section, we provide two propositions which will be crucial
 when evaluating the performance of our test procedure. They deal
with the behaviors of its expectation and its variance.

\begin{proposition}\label{prop2} Let $j$ be any given level parameter. Then,
{\small $$
\mathbb{E}_{_{\overrightarrow{p},\overrightarrow{q}}}(T_j)=\sum_{l=1}^{M}\sum_k\left( \int_\mathbb{R}(p_l-q_l)\phi_{jk}\right)^2-\frac{1}{n^2}
\sum_{l=1}^{M}\sum_k\sum_{i=1}^n\left(\int_{\mathbb{R}}\left(a_l(i)f_i-b_l(i)g_i\right)\phi_{jk}\right)^2.$$}
\end{proposition}

\begin{remark}\label{rem1}
For the particular case where the sequences of the mixing weights $(\omega_u(i), 1 \leq u \leq M, 1 \leq i \leq n)$ and $(\sigma_u(i), 1 \leq u \leq M, 1 \leq i \leq n)$
are identical, the test statistic $T_j$ is centered under the null hypothesis.
\end{remark}

\noindent 
\begin{corollary}\label{coro1}
For any $j \in \mathbb{N}$,
$$\left|\mathbb{E}_{_{\overrightarrow{p},\overrightarrow{q}}}(T_j)-\begin{displaystyle}\sum_{l=1}^{M}\sum_k\end{displaystyle}\left( \int_\mathbb{R}(p_l-q_l)\phi_{jk}\right)^2 \right| 
\leq \frac{8LMR^2}{Kn}.$$

\end{corollary}

\begin{proposition}\label{prop3}
\noindent There exists a constant $C_{_T}=C_{_T}(R,L,\|\phi\|_\infty)>0$ such that
\begin{eqnarray*}\label{varcal}
\mathbb{V}ar_{_{\overrightarrow{p},\overrightarrow{q}}}(T_j) \leq C_{_T}\left(\frac{2^j}{n^2}+\frac{1}{n}\sum_l\|p_l-q_l\|_2^2+ \sqrt{\frac{2^j}{n^3}}\sum_l\|p_l-q_l\|_2\right)\frac{M^2}{K^2}.
\end{eqnarray*}
\end{proposition}

\begin{remark}\label{rem2}
Under the null hypothesis the variance of the test statistic $T_j$ is
less than or equal to $C_{_T}M^2K^{-2} \ 2^j \ n^{-2}.$
\end{remark}

\subsection{Minimax performance of the test procedure}\label{minimaxna}

\noindent For any $s>0$, let $(r_n)_{n\in \mathbb{N}}$ be the sequence such that 
$$r_n=n^{-\frac{2s}{1+4s}}\quad \forall n\in \mathbb{N}^*.$$ 
The following theorem shows that the test procedure defined in 
section \ref{section2} provides an accurate upper bound when it is well
calibrated.

\begin{theorem}[Upper bound] \label{theoup}
Fix $0<\gamma<1$ and consider the test procedure $\Delta_s^*=\Delta_{j_n}$ where $j_n$ is the smallest integer such that $2^{-j_n} \leq n^{-\frac{2}{1+4s}}.$ Let $t$ and $C_\gamma$  be two positive real numbers defined as follows :
\begin{eqnarray*}
t=\left(2\sqrt{\frac{C_{_T}}{\gamma}} +8LR^2\right)\frac{M}{K},\\
C_\gamma^2=2\left(\frac{1}{K}\sqrt{\frac{6\ C_{_T}}{\gamma}}+R+\frac{t}{M}\right).
 \end{eqnarray*}
 \noindent  Then 
{\small{\begin{eqnarray}\label{uptheo}
\limsup_{n \to \infty}  \left( \sup_{(\overrightarrow{p},\overrightarrow{q}) \in \Theta_0(R)}\mathbb{P}_{\overrightarrow{p},\overrightarrow{q}}(\Delta_s^*=1) + \sup_{(\overrightarrow{p},\overrightarrow{q})  \in \Theta_1(R,C,n,s)} \mathbb{P}_{\overrightarrow{p},\overrightarrow{q}}(\Delta_s^*=0)\right) \leq \gamma
\end{eqnarray}}}
for all $C>C_\gamma$.
\end{theorem}

\noindent Although the exact value of the constant $C_T$ is very 
complicated, it can be exactly calculated by following the proofs.

\noindent Now, let us focus on the lower bound associated with our nonparametric testing problem $\mathcal{H}_0$ versus $\mathcal{H}_1$. \\
We aim at providing a constant $c_\gamma$ such that we ensure that
no test procedure is able to choose $\mathcal{H}_0$ or $\mathcal{H}_1$ with a sum of the probability errors less than $\gamma$ ($0<\gamma<1$). Obviously, the smaller the distance between $c_\gamma$ and $C_\gamma$  the more accurate our results.
The next theorem proves that our test procedure is asymptotically minimax.\\
Similarly to the classical methods for providing lower bounds (see for instance Gayraud and Pouet \cite{PG} or Butucea and Tribouley \cite{BT}) we shall consider a subspace of $\Lambda_n(R,C)$ that is,
for any chosen $C_1>0,$ 
\begin{equation}
\tilde{\Lambda}_n(R,C,C_1)=\left\{(p,q) \in \Lambda_n(R,C); \inf_{z \in [0,1[} \min (p(z), q(z))\geq C_1 \right\}.
\label{spaceinf}
\end{equation}

\begin{theorem}[Lower bound] \label{theolow}
Let $0<\gamma<1$, $s>0$  and let $c_\gamma>0$ satisfy  $$
c_\gamma^4 = \left( \frac {C_1^2} {L \ K^2} \ln[4(1-\gamma)^2+1] \wedge 2R^2\right)\frac{2^{-4s}}{4M^2}.
$$
Then   for all $C<c_\gamma$
{\small \begin{eqnarray}\label{lowtheo}
\liminf_{n \to \infty}  
\inf_{\Delta}
\left( \sup_{(\overrightarrow{p},\overrightarrow{q}) \in \Theta_0(R)}\mathbb{P}_{\overrightarrow{p},\overrightarrow{q}}(\Delta=1) + \sup_{(\overrightarrow{p},\overrightarrow{q})  \in \Theta_1(R,C,n,s)} \mathbb{P}_{\overrightarrow{p},\overrightarrow{q}}(\Delta=0)\right) > \gamma
\end{eqnarray}}
where the infimum is taken over all test procedure $\Delta$.
\end{theorem}

\noindent  From Theorems \ref{theoup} and \ref{theolow} we deduce
the minimax rate of testing. It is the same as the one found by 
Butucea and Tribouley \cite{BT} when there is only one subgroup.
Advances in our results 
are the extension to the varying mixing weights model which 
allows non-identically distributed random variables compared
to Butucea and Tribouley \cite{BT} and
 the role played by the mixing weights
which is clearly exposed.

\begin{corollary} For any $s>0$, the test procedure $\Delta_s^*$ is asymptically minimax and the minimax rate separating $\mathcal{H}_0$ and $\mathcal{H}_1$ is $r_n=n^{-\frac{2s}{1+2s}}$. 
\end{corollary}


\subsection{Discussion about the constants $c_\gamma$ and $C_\gamma$}

\noindent In the two previous theorems we exhibited two constants appearing in
 the upper and the lower bounds. We think that the connection between these constants and the model's parameters $M$ and $K$ is a novelty and really deserves a discussion. Indeed, we keep in mind that 
\begin{itemize}
\item $C_\gamma$ is the minimal value for $C$ such that
 our test statistic is able to detect if all the mixture components are identical in the two populations with the sum of the probability errors not exceeding $\gamma$;
\item $c_\gamma$ is the maximal value for $C$ such that no test statistic is able to detect if all the mixture components are identical in the two populations with the sum of probability errors not exceeding $\gamma$.
\end{itemize} 

\noindent As a consequence we proved that our test statistic is optimal in the minimax sense since it attains the minimax rate of convergence separating $\mathcal{H}_0$ and $\mathcal{H}_1$.\\

\noindent According to the definitions of $c_\gamma$ and $C_\gamma$
we let the reader be aware that:
\begin{itemize}
\item the smaller the constant $K$, the larger the family of the
mixing weights satisfying HYP-1;
\item the smaller the constant $M$, the bigger ($=$ the worse)  the constant $C_\gamma$ and the bigger  the constant $c_\gamma$;
\item the smaller the constant $K$, the bigger ($=$ the worse)  the constant $C_\gamma$ and the bigger  the constant $c_\gamma$.
\end{itemize}

\noindent Although the exact separation constant is not 
established in this study (since $c_\gamma \not= C_\gamma$), we prove that $c_\gamma$ and $C_\gamma$ strongly depend on the  smallest eigenvalue of the matrices $\Omega\Omega^*$ and $\Sigma\Sigma^*$. 

\section{Numerical experiments and application}\label{section3}

The aim of this section is twofold: to illustrate by numerical experiments the good performance of the test procedures based on the statistics $T_{j_n}$ and to show the usefulness of our method on real data.\\

\noindent First,  $2$ examples of mixture models are given 
to show the interest of
the problem we have considered. Next we illustrate the behaviour of the
test statistics $T_{j_n}$.

\subsection{Examples of mixture models}
\noindent {\underline{Figure \ref{Fig2M}}: [Mixture with two components]\\
\noindent Consider two populations sampled from the same mixture densities such that 
\begin{itemize}
\item the size of the two populations $(Y,Z)$ is $n=500$,
\item the ranks of the matrices of the mixing weights $\Omega^*$ and $\Sigma^*$ are $2$,
\item the two components of the mixtures are the uniform density 
$\mathcal{U}\left([-1,0]\right)$ and the normal density $\mathcal{N}(3,4).$ 
\end{itemize}
\noindent 

\setcounter{figure}{0}
\begin{figure}[h]
\begin{center}
\resizebox*{6cm}{!}{\includegraphics{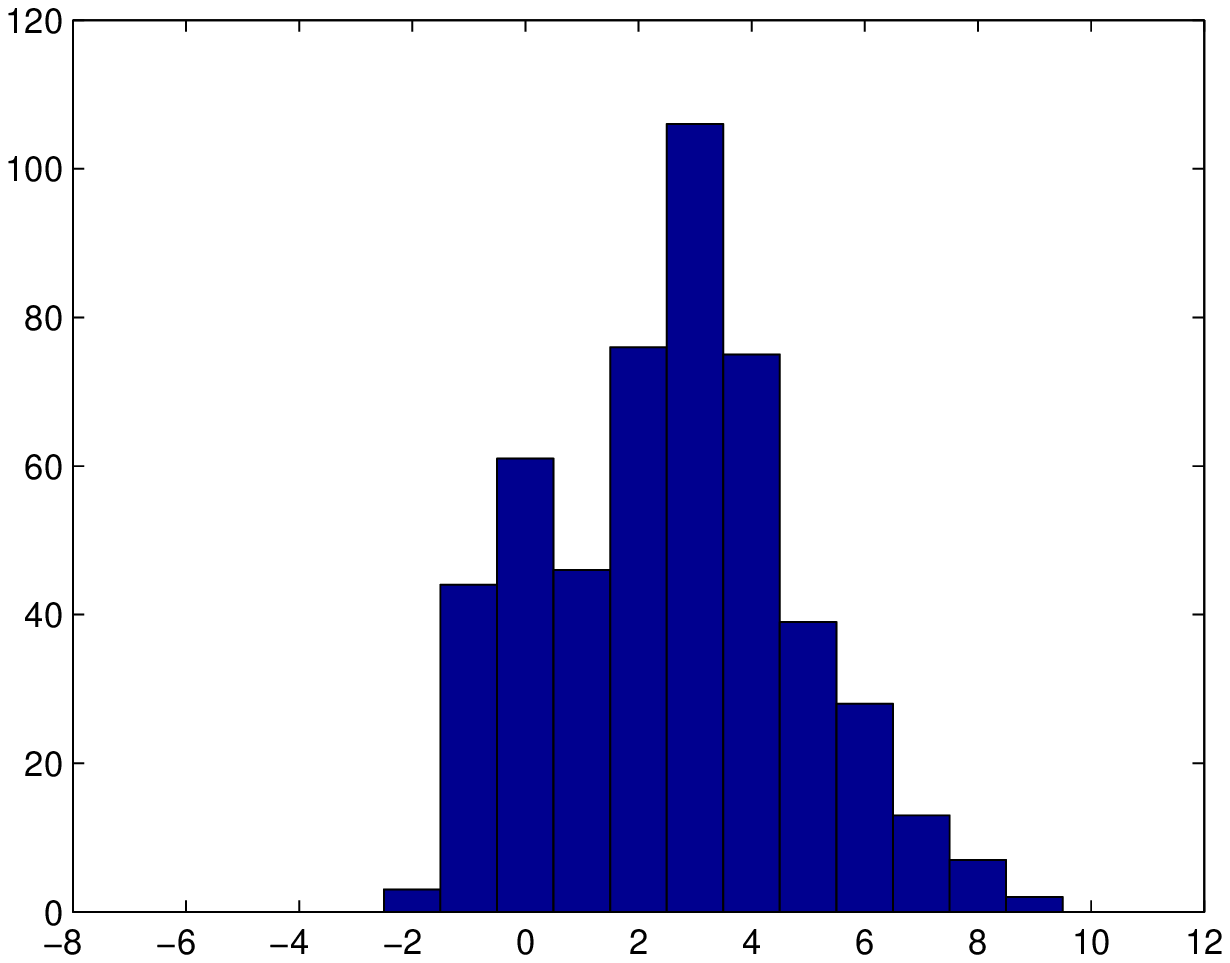}}
\resizebox*{6cm}{!}{\includegraphics{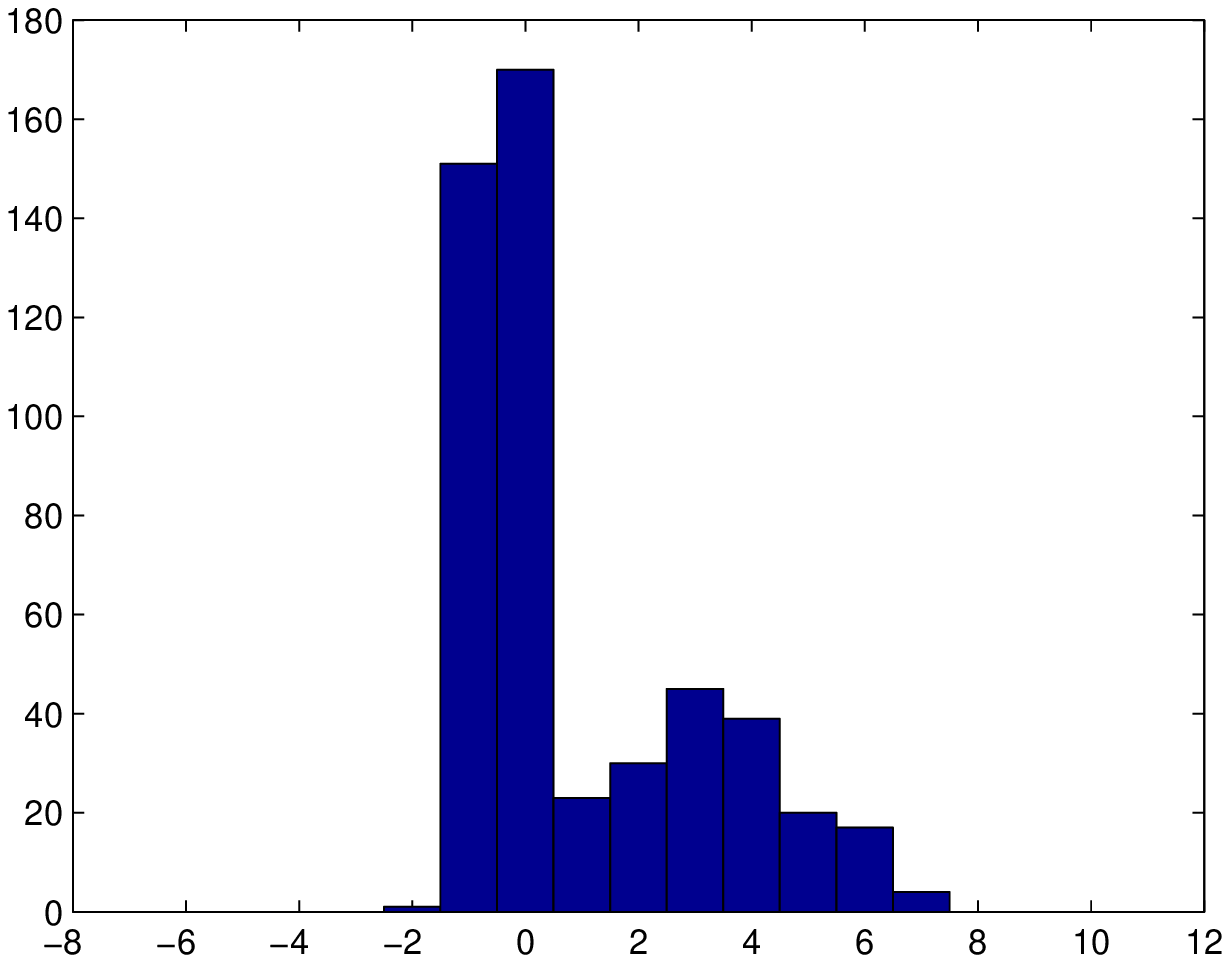}}
\caption{Histogram (a) of population $Y$ and histogram (b) population $Z$.}
\label{Fig2M}
\end{center}
\end{figure}

\noindent {\underline{Figure \ref{Fig3M}}}: [Mixture with three components]\\
\noindent Consider two populations sampled from the same mixture densities such that 
\begin{itemize}
\item the size of the two populations ($Y,Z$) is $n=500$,
\item the ranks of the matrices of the mixing weights $\Omega^*$ and $\Sigma^*$ are $3$,
\item the three components of the mixtures are the normal densities
$\mathcal{N}(-2,1),$ $\mathcal{N}(0,1)$ and $\mathcal{N}(2,1)$. 
\end{itemize}
\noindent

\setcounter{figure}{1}
\begin{figure} [h]
\begin{center}
\resizebox*{6cm}{!}{\includegraphics{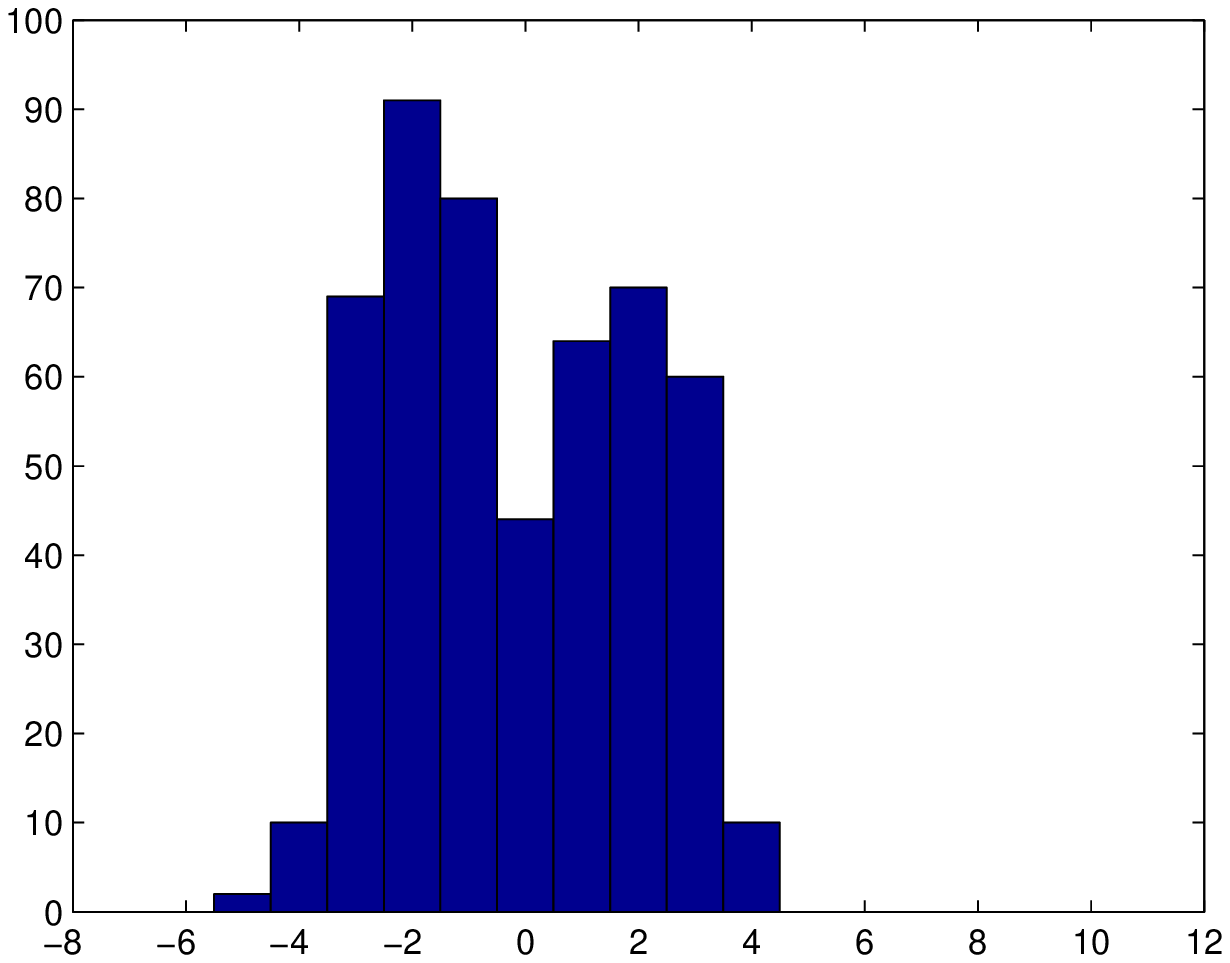}}
\resizebox*{6cm}{!}{\includegraphics{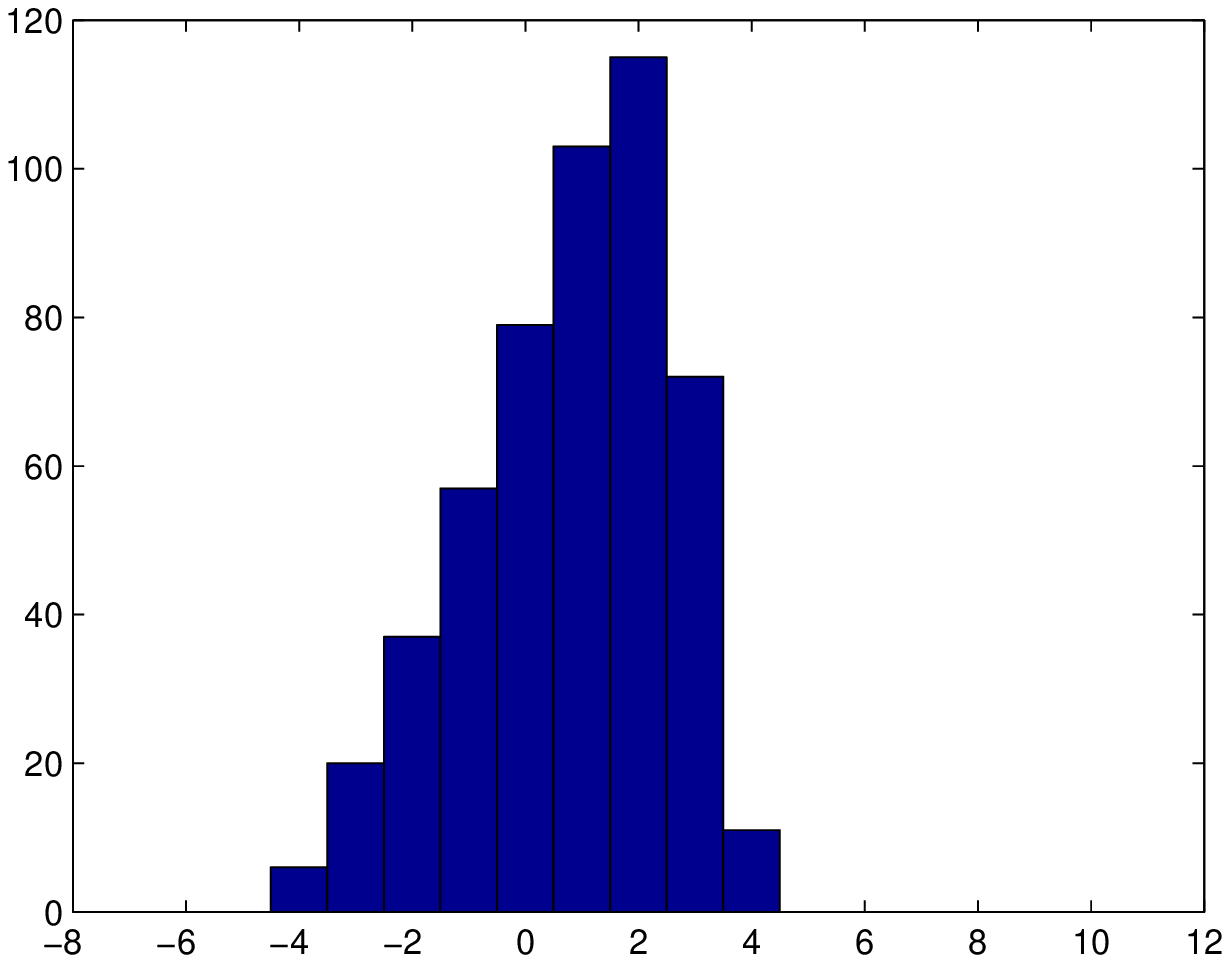}}
\caption{Histogram (a) of population $Y$ and histogram (b) of population $Z$.}
\label{Fig3M}
\end{center}
\end{figure}

\newpage
\noindent The histograms of the observations are quite different in Figures $1$ and $2$, although they correspond to mixture models with the same components. So the previous schemes show how hard it is to guess whether the mixture components of the two populations $(Y,Z)$ are exactly the same or not. Hence, it justifies that the statistician needs an adequate test statistic to decide whether the populations $(Y,Z)$ have the same mixture components or not.

\subsection{Construction of the test procedure: calibration of $t_n$}

\noindent In the theoretical part of this paper we provide a decision rule to test $\mathcal{H}_0$ against $\mathcal{H}_1$. This decision rule $\Delta_{j_n}$ relies on the sign of  $T_{j_n}-t_n,$ where $t_n$ is the threshold value depending on the sum of the errors $\gamma$ and $T_{j_n}$ is the test statistic. In the positive case (resp. in the negative case) $\Delta_{j_n}$ proposes to accept $\mathcal{H}_1$ (resp. $\mathcal{H}_0$).  

\noindent From the practical point of view, we give some hints
to adjust the threshold value $t_n$. Here we use the Haar basis and we set $s=4$. For this, we consider two different approaches.\\

\noindent The {\it{first}} approach consists in fixing the first type error, 
$0<\gamma_1<1$, and in choosing $t_n$ as the quantile of order $1-\gamma_1$
of the test statistic obtained after $1000$ replications of the chosen 
mixture model.\\

\noindent The {\it{second}} approach consists in choosing $t_n$ as the value for which the sum of the two errors is the minimal one according 
of the statistic of test obtained after $1000$ replications of the mixture model chosen.


\subsection{Connection between $K$ and the performance of the test procedure.}

\noindent The aim of this paragraph is to illustrate the connection between the value of $K$ and the performance of our test procedure. We provide simulations of  Gaussian mixture models and we give for  several values of $n$
\begin{itemize}
\item the value of $t_n$ associated with a first type error equal to $10\%$,
\item the power of the test procedure based on the threshold value $t_n$,
\item the minimum of the global error $\gamma_{opt}$ - the sum of the first type and the second type errors - reachable by the test procedure,
\item the value $t_{opt}$ which corresponds to the global error $\gamma_{opt}$.
\end{itemize}

\noindent 
We consider two samples: $Y_1,\ldots,Y_n$ and $Z_1,\ldots, Z_n$.
Two mixture components are such that
\begin{itemize}
\item under $H_0$, $p_1(\cdot) = q_1(\cdot)
 \sim \mathcal{N}(-2,1)$ and $p_2(\cdot) =q_2(\cdot) \sim \mathcal{N}(3,4)$,
\item under $H_1$,
$p_1(\cdot)  \sim \mathcal{N}(-2,1)$, $p_2(\cdot) \sim \mathcal{N}(3,4)$,
 $q_1(\cdot) \sim \mathcal{N}(0,1)$ and $q_2(\cdot) \sim \mathcal{N}(1,1)$.
\end{itemize}

\noindent Weights of samples $Y$ and $Z$ for
 Gaussian Model $1$ are described in
 Table $1$. \\
\begin{center}
\begin{tabular}{|c|l|c|c|} \hline
  Sample & Range of $i$ & $\sigma_1(i)$ or $\omega_1(i)$ &
$\sigma_2(i)$ or $\omega_2(i)$ \\ \hline
Y & $i=1,\ldots,0.8\ n$ & $0.6$ & $0.4$ \\ \cline{2-4}
  & $i=0.8\ n +1,\ldots,n$ & $0.4$ & $0.6$ \\ \hline
Z & $i=1,\ldots,0.3\ n$ & $0.2$ & $0.8$ \\ \cline{2-4}
  & $i=0.3\ n+1,\ldots,\ n$ & $0.5$ & $0.5$ \\ \hline
\end{tabular}
\end{center}
\begin{center}
Table $1$: Model $1$
\end{center}

\noindent The results are given in  Table $2$. We point out that the
constant $K$ related to the smallest eigenvalue is very close to $0$.
Therefore we expect poor results.\\
\begin{center}
\begin{tabular}{|c|c|c|c|}
\hline  Gaussian Model $1$ &$n=200$&$n=500$&$n=1000$\\
\hline 
 $t_n$&0.289&0.135&0.080\\
$\hbox{Power}$&36.7\%&68.1\%&85.7\%\\
\hline
$\gamma_{opt}$&52.6\%&38.2\%&23.5\%\\
$t_{opt}$&0.022&0.080&0.092\\
\hline
\end{tabular}
\end{center}
\begin{center}
Table $2$: $K=0.013$\\
\end{center}

\noindent 
 Weights of samples $Y$ and $Z$ for
 Gaussian Model $2$ are described in
 Table $3$. \\
\begin{center}
\begin{tabular}{|c|l|c|c|} \hline
  Sample & Range of $i$ & $\sigma_1(i)$ or $\omega_1(i)$ &
$\sigma_2(i)$ or $\omega_2(i)$ \\ \hline
Y & $i=1,\ldots,0.8\ n$ & $0.8$ & $0.2$ \\ \cline{2-4}
  & $i=0.8\ n +1,\ldots,n$ & $0.3$ & $0.7$ \\ \hline
Z & $i=1,\ldots,0.3\ n$ & $0.1$ & $0.9$ \\ \cline{2-4}
  & $i=0.3\ n+1,\ldots,\ n$ & $0.4$ & $0.6$ \\ \hline
\end{tabular}
\end{center}
\begin{center}
Table $3$: Model $2$
\end{center}

\noindent For this setup, the constant $K$ is almost three times the one appearing
in Gaussian Model 1. Therefore we expect improved results.\\
\begin{center}
\begin{tabular}{|c|c|c|c|}
\hline  Gaussian Model $2$ &$n=200$&$n=500$&$n=1000$\\
\hline 
 $t_n$&0.994&0.061&0.027\\
$\hbox{Power}$&85.2\%&91.5\%&96.8\%\\
\hline
$\gamma_{opt}$&24.6\%&16.3\%&9.5\%\\
$t_{opt}$&0.078&0.103&0.047\\
\hline 
\end{tabular}
\end{center}
\begin{center}
Table $4$: $K=0.033$
\end{center}

\noindent  Weights of samples $Y$ and $Z$ for
 Gaussian Model $3$ are described in
Table $5$. \\
\begin{center}
\begin{tabular}{|c|l|c|c|} \hline
  Sample & Range of $i$ & $\sigma_1(i)$ or $\omega_1(i)$ &
$\sigma_2(i)$ or $\omega_2(i)$ \\ \hline
Y & $i=1,\ldots,0.8\ n$ & $0.8$ & $0.2$ \\ \cline{2-4}
  & $i=0.8\ n +1,\ldots,n$ & $0.3$ & $0.7$ \\ \hline
Z & $i=1,\ldots,0.3\ n$ & $0.9$ & $0.1$ \\ \cline{2-4}
  & $i=0.3\ n+1,\ldots,\ n$ & $0.3$ & $0.7$ \\ \hline
\end{tabular}
\end{center}
\begin{center}
Table $5$: Model $3$
\end{center}

\vspace{0.5cm}
\noindent In this setup, the constant $K$ is more than five times the one
appearing in Gaussian Model 1 and more than twice the one
appearing in Gaussian Model 2. Therefore we expect better results.\\
\begin{center}
\begin{tabular}{|c|c|c|c|}
\hline  Gaussian Model $3$ &$n=200$&$n=500$&$n=1000$\\
\hline 
 $t_n$&0.054&0.030&0.015\\
$\mathcal{P}$&97.1\%&96.7\%&98.1\%\\
\hline
$\gamma_{opt}$&10.5\%&9.6\%&6.5\%\\
$t_{opt}$&0.066&0.064&0.034\\
\hline
\end{tabular}
\end{center}
\begin{center}
Table $6$ : $K=0.068$
\end{center}

\noindent According to numerical results in Tables $2,4$ and $6$, it is clear that for a fixed $n$, the larger the value of $K$, the better the performance of the test procedure. Indeed, when the first type error is $10\%$, we see that  increasing values of $K$ increases the power of the test procedure.
 Moreover, we remark that 
the optimal global error $\gamma_{opt}$ increases when the value of $K$ decreases. In fact, this is not surprising as this behaviour
was predicted by our theoretical results: the smaller the value of $K$ the larger the constant $C_\gamma$ (see Theorem \ref{theoup}). In other words, in a mixture model with a small value of $K$ one needs a lot of observations to ensure  good performance of our test procedure. 

\subsection{Application to real data}

\noindent In this part we apply our results to real data. The dataset comes from 
a survey conducted by the french national statistical agency called InstituT National de Statistique et d'Etudes Economiques (abbreviated to INSEE). This survey 
called {\it{D\'eclaration Annuelle des Donn\'ees Sociales}} (abbreviated to DADS) took place in 2007 and is about employees and related variables such as salary, working time or type of jobs. All information regarding this survey
can be found on the website of INSEE (see DADS 2007 postes et salari\'es, http://www.insee.fr.).
As far as we are concerned, we focused on working time per year. More precisely our goal is to make two
comparisons at the same time:
\begin{enumerate}
\item working time of men in Ile-de-France (region surrounding
Paris in France, abbreviated to $\mathcal{I}$ below) and the one done by men in all other regions of France (abbreviated to $\mathcal{P}$ below),
\item working time of women in Ile-de-France and the one done by women in all other regions of France.
\end{enumerate}
In this study we decide to only consider highly skilled workers such as
executive staff, managers. There are two populations:
\begin{itemize}
\item commercial and administrative staff (abbreviated to {\it{CAd}}),
\item technical staff (abbreviated to {\it{Tech}}).
\end{itemize}
We restrict to people working more than $1 \ 645$ hours per year.
The variable of interest is the number of working hours per year
divided by $1 \ 645$. Therefore it is a ratio equals to or greater than $1$.\\
 Available information
about different subpopulations of $\mathcal{I}$ and  $\mathcal{P}$ is gathered in the following table:

\begin{center}
\begin{tabular}{|c|c|c|c|c|}
\hline
  &    \multicolumn{2}{|c}{Ile-de-France ($\mathcal{I})$} & \multicolumn{2}{|c|}{Other regions ($\mathcal{P})$}  \\
   \hline
 Executive staff & Men & Women & Men & Women \\ \hline \hline
 \it{CAd} & $58.99\%$ & $41.01\%$ & $67.96\%$ & $32.04\%$ \\ \hline
 \it{Tech} & $81.08\%$ & $18.92\%$ & $86.72\%$ & $13.28\%$ \\ \hline
\end{tabular}
\end{center}
\begin{center}
Table $7$: Proportions of subpopulations by sex, area and job\\
\end{center}

\noindent 
There are $65\ 558$ people in $\mathcal{I}$ and  $75\ 062$ people in $\mathcal{P}$.\\

\noindent To begin, we pay attention to  the {\it{mean of the working-ratio}} of each population, namely $m_{\mathcal{I}}$ and $m_{\mathcal{P}}$.  Although information about sex (men or women) is available in the study conducted by INSEE, 
we assume that it 
is unknown in order to show the interest of our model. \\ 

\noindent Let $\sigma_{\mathcal{I}}$ and $\sigma_{\mathcal{P}}$ denote 
the standard deviations of population $\mathcal{I}$ and $\mathcal{P}$ according to the variable of interest. 
We suppose that a random sampling of order $n=5 \ 000$ in each population is available and is conducted as follows:
\begin{itemize}
\item $2 \ 500$ people living in $\mathcal{I}$  are  {\it{CAd}}  and $2 \ 500$ people living in $\mathcal{I}$ are {\it{Tech}}, 
\item $2 \ 500$ people living in $\mathcal{P}$  are  {\it{CAd}}  and $2 \ 500$ people living in $\mathcal{P}$ are {\it{Tech}} .
\end{itemize}
We are interested in the preliminary testing problem ($\mathcal{T}_{1}$):
$$ \mathcal{H}_{0} : m_{\mathcal{I}}=m_{\mathcal{P}} \quad vs \quad  \mathcal{H}_{1} : m_{\mathcal{I}} \not= m_{\mathcal{P}}.$$

\noindent We decide to address this testing problem by using the test statistic $$U=\frac{|\hat{m}_{\mathcal{I}}-\hat{m}_{\mathcal{P}}|}{\sqrt{\hat{\sigma}_{\mathcal{I}}^2+\hat{\sigma}_{\mathcal{P}}^2}},$$ 

\noindent where $\hat{m}_{\mathcal{I}}$ (resp. $\hat{m}_{\mathcal{P}}$) and $\hat{\sigma}_{\mathcal{I}}$ (resp. $\hat{\sigma}_{\mathcal{P}}$)
denote the usual estimators of $m_{\mathcal{I}}$ (resp. $m_{\mathcal{P}}$) and $\sigma_{\mathcal{I}}$ (resp. $\sigma_{\mathcal{P}}$), when using stratified random samplings like ours. Under the null hypothesis $\mathcal{H}_{0}$, the random variable $U$ is asymptotically normally distributed with mean $0$  and  variance $1$. \\

\noindent Here are the values computed from the samples: 

\begin{center}
\begin{tabular}{|c|c|}
\hline
 Ile-de-France ($\mathcal{I}$)& Other regions ($\mathcal{P}$)\\
\hline 
$\hat{m}_{\mathcal{I}}=1.1605$&$\hat{m}_{\mathcal{P}}=1.1531$\\
$\hat{\sigma}_{\mathcal{I}}=0.0015$&$\hat{\sigma}_{\mathcal{P}}=0.0014$\\
\hline 
\end{tabular}
\end{center}
\begin{center}
Table $8$: Estimated means and standard deviations by area
\end{center}

\noindent  The value of the test statistic $U$ is $3.5582$.  
The related $p$-value is close to $0.0026$. 
According to that, it strongly seems that $m_{\mathcal{I}} \not=m_{\mathcal{P}}$.  In other words,  $\mathcal{H}_{0}$ is rejected. \\\\

\noindent At this stage, a natural question arises : what is the reason of such a difference? Two hypotheses could explain  it:
\begin{enumerate}
\item distincts values of $m_{\mathcal{I}}$ and $m_{\mathcal{P}}$  are only related to the different proportions of men (or analogously women) between the two populations:
\begin{center}
\begin{tabular}{|c|c|c|}
\hline
 & Ile-de-France ($\mathcal{I}$)& Other regions ($\mathcal{P}$) \\
   \hline
Men & $ 68.93\% \ (45187)$ &  $ 76.70\% \ (57575)$ \\ \hline
Women & $31.07\% \ (20371)$ & $23.30\% \ (17487)$ \\ \hline
\end{tabular}
\end{center}
\begin{center}
Table $9$: Proportions of subpopulations by area and sex
\end{center}
\item distincts values of $m_{\mathcal{I}}$ and $m_{\mathcal{P}}$ are  also related to  different distributions of  working-ratio  of population $\mathcal{I}$
(abbreviated to $W.R.^{(\mathcal{I})}$) and working-ratio 
 of  population $\mathcal{P}$ (abbreviated to $W.R.^{(\mathcal{P})}$).
\end{enumerate}


\noindent Trusting one of these new hypotheses becomes at first glance difficult to argue when only considering two random samples of size $n$ in each population without the knowledge of sex (man or woman). Nevertheless, a way to address the testing problem ($\mathcal{T}_{2}$):
\begin{eqnarray*}&\mathcal{H'}_{0} :& \hbox{distributions of }  W.R.^{(\mathcal{I})} \hbox{and }  W.R.^{(\mathcal{P})}  \hbox{conditionnally to sex are identical} \\
vs \quad &\mathcal{H'}_{1} :& \hbox{distributions of } W.R.^{(\mathcal{I})} \hbox{and }  W.R.^{(\mathcal{P})}  \hbox{conditionnally to sex are different}
\end{eqnarray*}
\noindent is to consider our testing procedure. \\

\noindent Let $p_{1}$ and $p_{2}$ (resp. $q_{1}$ and $q_{2}$) denote the density functions of the random variables $W.R.^{(\mathcal{I})}|_{man}$ and $W.R.^{(\mathcal{I})}|_{woman}$ (resp. $W.R.^{(\mathcal{P})}|_{man}$ and $W.R.^{(\mathcal{P})}|_{woman}$).\\

\noindent The testing problem $\mathcal{T}_{2}$ can be written as follows:
$$ \mathcal{H'}_{0} : p_{1}=q_{1}  \hbox{ and } \ p_{2}=q_{2} \quad vs \quad  \mathcal{H'}_{1} : p_{1}\not=q_{1}  \hbox{ or } \ p_{2}\not=q_{2}.$$

\noindent Observations of the working-ratio random variables $Y_{1}, \dots, Y_{n}$ (resp. $Z_{1}, \dots, Z_{n}$) in 
population $\mathcal{I}$ (resp. in $\mathcal{P}$) are available. The mixture model we get is the one described in Section 
$1.3$ with:
\begin{itemize}
\item $M=2$ and $n=5 \ 000$,
\item $\left(\omega_{1}(i),\omega_{2}(i)\right)=(0.5899,0.4101)$  for a $n/2$-tuple of indices,
\item $\left(\omega_{1}(i),\omega_{2}(i)\right)=(0.8108,0.1892)$  for a $n/2$-tuple of indices,
\item $\left(\sigma_{1}(i),\sigma_{2}(i)\right)=(0.6796,0.3204)$  for a $n/2$-tuple of indices,
\item $\left(\sigma_{1}(i),\sigma_{2}(i)\right)=(0.8672,0.1328)$  for a $n/2$-tuple of indices.
\end{itemize}

\noindent Let us describe the methodology of the testing procedure applied 
to these real data. We use the test
studied in Section $2$ with regularity parameter
$s=4$ and choose the usual Haar wavelet to construct our test statistic $T_{j_{s}}$. 
The threshold value of the testing
procedure is computed according to the following heuristics:
$ t = s t_\alpha$ where $t_\alpha$ is the $1-\alpha$ Gaussian quantile and
$s$ is the standard deviation of the test statistics estimated by bootstrap (resampling is made $200$ times).
As we choose $\alpha=10\%$, we have $t_{0.1} = 1.28$.\\

\noindent  The value of  $T_{j_{s}}$ obtained is 
$t_{j_{s}}=0.5412$ whereas the threshold value is $t=0.3324$. Since $t_{j_{s}}$ is larger than the threshold value $t$, we conclude that there exists a difference between the distributions  $W.R.^{(\mathcal{I})}$ and $W.R.^{(\mathcal{P})}$ conditionnally to sex. In other words,  $\mathcal{H'}_{0}$ is rejected. \\

\noindent In this last paragraph, we study the numerical performances of our testing procedure, built from $T_{j_{s}}$.  For several values of $n$, a sample of size $n$ is drawn from $\mathcal{I}$ (resp. $\mathcal{P}$) and is divided into two subsamples :
one subsample  of size $n/2$ is drawn from the subpopulation {\it{CAd}} and
the other is drawn from the subpopulation {\it{Tech}}. \\

\noindent For each value of $n$, $200$ samples are drawn.
 The results are gathered in the following table: 

\begin{center}
\begin{tabular}{|c|c|c|c|}
\hline 
Sample size n & First type error: $E_{\mathcal{I}}^{(n)}$ & First type error: $E_{\mathcal{P}}^{(n)}$
 & Power \\
 \hline 
1\ 000 &  0 & 0 & $0.110$\\
2\ 000 & 0 & $0.005$ & $0.185$ \\
3\ 000 & $0.005$ & $0.005$ & $0.335$ \\
4\ 000 & 0 & $0.005$ & $0.530$ \\
\hline
5\ 000 & $0.005$ & 0 & $0.635$  \\
6\ 000 & $0.005$ & 0 & $0.745$ \\
8\ 000 & $0.005$ & $0.005$ & $0.925$\\
\hline
\end{tabular}
\end{center}
\begin{center}
Table $10$: First type error and power of the method
\end{center}

\noindent 
$-$ First type error $E_{\mathcal{I}}^{(n)}$ is  the proportion of observations of $T_{j_{s}}$ larger than the threshold value, when comparing two samples of size $n$ in $\mathcal{I}$.\\
$-$ First type error $E_{\mathcal{P}}^{(n)}$ is the proportion of observations of $T_{j_{s}}$ larger than the threshold value, when comparing two samples of size $n$ in $\mathcal{P}$.\\
$-$ Power is the proportion of observations $T_{j_{s}}$ larger than the threshold value, when comparing a sample of size $n$ in $\mathcal{I}$ and a sample of size $n$ in $\mathcal{P}.$\\ 

\noindent 
It appears that the testing procedure with the heuristically chosen threshold is very
conservative. This is the only drawback of our methodology. Nevertheless the
behaviour of the testing procedure is as expected: the larger the sample
size the larger the power. As we see, for the cases $n \geq 5 \ 000$, 
our testing procedure is powerful. It tends to prove that there exists a difference between  the working-ratios of the two populations conditionally to sex. \\

\noindent This study on DADS 2007 demonstrates the usefulness of the
varying mixing weights model. It really suggests
that our testing procedure can be successfully applied to all types of 
data in social science. From our point of view, researchers in social science 
should consider the mixing varying weights model and our testing procedure 
as soon as some information at the individual level has been omitted during
 a survey and is available at higher levels.

\section{Open questions}\label{section4}

\noindent As a conclusion, we have provided a statistical procedure for a testing problem on the mixture components of  two populations $\left( Y, Z \right)$. This one was proved to be optimal in the minimax sense (Theorems \ref{theoup} and \ref{theolow}). In addition, we explained how the weights of the mixture model influence the performance of the statistical rule. All these theoretical results are illustrated by our numerical experiments.    \\

\noindent It seems to us important to give some hints about possible extensions of this work. From the theoretical and practical
points of view, it would be interesting to study the same problem without assuming that the mixing weights are exactly known to the statistician.
Several explanations can be given
\begin{itemize}
\item the statistician can estimate the mixing weights for 
an observation by using covariates and an appropriate predictive model
such as the logistic one,
\item a Bayesian approach is chosen for
the mixing weights,
\item exogenous information allows the statistician
to roughly estimate the mixing weights.
\end{itemize}
In this case several natural questions arise
\begin{itemize}
\item What statistical rule should be considered?
\item What kind of performance can be expected for such a rule?
\item How much do random mixing weights deteriorate the performance? 
\end{itemize}
Such questions are beyond the scope of this article and their answers
certainly involve random matrices theory. \\

\noindent Finally, it would be nice to show how to choose the 
adequate value of $t_n$ in a better way than the complicated one given in 
Theorem \ref{theolow}.

\section{Proofs of main results}\label{section5}

\noindent This section is devoted to the proofs of our results. 
The proofs often need technical lemmas which shall be proved in Appendix. For the sake of simplicity we sometimes omit
$\overrightarrow{p}$ and $\overrightarrow{q}$
in the indices when there is no ambiguity.

\subsection{Proofs of Propositions and Corollaries}
\noindent {\it{Proof of Proposition \ref{maibo}}}: We refer to Maiboroda \cite{Mai-3}. A solution of the two problems is given for any  $(l,i) \in  \{1, \dots, M\} \times \{1, \dots, n\}$ by $$a_l(i)= \frac{1}{det(\Gamma_n)}\sum_{u=1}^M(-1)^{l+u}\gamma_{_{lu}}\omega_u(i)$$
$$b_l(i)= \frac{1}{det(\Gamma'_n)}\sum_{u=1}^M(-1)^{l+u}\gamma'_{_{lu}}\sigma_u(i)$$
where $\gamma_{_{lu}}$ and $\gamma'_{_{lu}}$ are respectively the minor $(l,u)$ of the matrix $\Gamma_n$ and  the minor $(l,u)$ of the matrix $\Gamma'_n$. 
Inequalities (\ref{hypflo1}) and (\ref{hypflo2}) are obtained by using lemma \ref{lem0}.\\

\noindent 
\hfill
$\Box$ \\

\noindent {\it{Proof of Proposition \ref{prop2}}}: 
Let us evaluate the expectation of $T_j$. 
{\small{\begin{eqnarray*}
&&\mathbb{E}_{\overrightarrow{p},\overrightarrow{q}}(T_j) =  \mathbb{E}_{\overrightarrow{p},\overrightarrow{q}}\left(
\frac 1 {n^2} \sum_{l=1}^M \sum_k \sum_{i_1 \neq i_2}
 (a_l(i_1)\phi_{jk}(Y_{i_1}) -b_l(i_1) \phi_{jk}(Z_{i_1}))
(a_l(i_2)\phi_{jk}(Y_{i_2}) -b_l(i_2) \phi_{jk}(Z_{i_2})) \right) \\
 & = & \frac 1 {n^2}  \sum_{l=1}^M \sum_{k} \sum_{i_1 \neq i_2}
\mathbb{E}_{\overrightarrow{p},\overrightarrow{q}}\left[ a_l(i_1)\phi_{jk}(Y_{i_1}) -b_l(i_1) \phi_{jk}(Z_{i_1}) \right]
\mathbb{E}_{\overrightarrow{p},\overrightarrow{q}}\left[ a_l(i_2)\phi_{jk}(Y_{i_2}) -b_l(i_2) \phi_{jk}(Z_{i_2}) \right],
\end{eqnarray*}}}
since the random variables $(Y_{i_1}, Z_{i_1})$ and $(Y_{i_2}, Z_{i_2})$
are independent.

\noindent We have for all $1 \leq i \leq n$,
$$
\mathbb{E}_{\overrightarrow{p},\overrightarrow{q}}\left[ a_l(i)\phi_{jk}(Y_{i}) - b_l(i)\phi_{jk}(Z_{i}) \right] 
 = \int_{\mathbb{R}} \left( \sum_{u=1}^M
\left(a_l(i)\omega_u(i)  p_u - b_l(i)\sigma_u(i) q_u \right) \right)\phi_{jk}.
$$

\noindent 
By introducing the diagonal term $i_1 = i_2$ in the sum, we get
\begin{eqnarray*}
\mathbb{E}_{\overrightarrow{p},\overrightarrow{q}}(T_j) & = & \frac 1 {n^2}  \sum_{l=1}^M \sum_k
\left(  \int_{\mathbb{R}} \phi_{jk} \left( \sum_{i=1}^n \sum_{u=1}^M
a_l(i)\omega_u(i) p_u - \sum_{i=1}^n \sum_{u=1}^Mb_l(i)\sigma_u(i)q_u  \right) \right)^2 \\
 & & \quad \quad \quad \quad \quad \quad \quad \quad \mathop{-}\frac{1}{n^2}
\sum_{l=1}^{M}\sum_k\sum_{i=1}^n\left(\int_{\mathbb{R}}\left(a_l(i)f_i-b_l(i)g_i\right)\phi_{jk}\right)^2\\
& = & \sum_{l=1}^{M}\sum_k\left( \int_\mathbb{R}
(p_l-q_l)\phi_{jk}  \right)^2 
 \mathop{-} \frac{1}{n^2}
\sum_{l=1}^{M}\sum_k\sum_{i=1}^n\left(\int_{\mathbb{R}}\left(a_l(i)f_i-b_l(i)g_i\right)\phi_{jk}\right)^2,
\end{eqnarray*}
because of the two properties
 $\begin{displaystyle}\frac 1 n \sum_{i=1}^n \end{displaystyle}a_l(i) \omega_u(i) = \delta_{lu}$ and $\begin{displaystyle}\frac 1 n \sum_{i=1}^n \end{displaystyle}b_l(i) \sigma_u(i) = \delta_{lu}$. Thus the result for the expectation is proved.\hfill
$\Box$ \\

\noindent {\it{Proof of Corollary \ref{coro1}}}:\\ 
 According to  proposition \ref{prop2} we only have to bound
the quantity $$
\frac{1}{n^2}
\sum_{l=1}^{M}\sum_k\sum_{i=1}^n\left(\int_{\mathbb{R}}\left(a_l(i)f_i-b_l(i)g_i\right)\phi_{jk}\right)^2.$$ 

\noindent Using the Cauchy-Schwarz inequality and lemma \ref{lem2}, we have

\begin{eqnarray*}
\sum_{l=1}^{M}\sum_k\sum_{i=1}^n\left(\int_{\mathbb{R}}\left(a_l(i)f_i-b_l(i)g_i\right)\phi_{jk}\right)^2 &\leq& 
\sum_{l=1}^{M}\sum_k\sum_{i=1}^n\int_{I_{jk}}\left(a_l(i)f_i-b_l(i)g_i\right)^2\int\phi^2_{jk}\\
&=&\sum_{l=1}^{M}\sum_{i=1}^n\left[\sum_k\int_{I_{jk}}\left(a_l(i)f_i-b_l(i)g_i\right)^2\right]\\
&\leq&2\sum_{i=1}^n\sum_{l=1}^{M}\left[\sum_k\int_{I_{jk}}\left(a_l(i)f_i\right)^2+ \int_{I_{jk}}\left(b_l(i)g_i\right)^2\right]\\
&\leq& 4L\left(\sum_{i=1}^{n}\sum_{l=1}^M a^2_l(i)\|f_i\|_2^2+ \sum_{i=1}^{n}\sum_{l=1}^Mb^2_l(i)\|g_i\|_2^2\right) \\
&\leq&\frac{8LMR^2n}{K}.
\end{eqnarray*}

\noindent Last inequality is due to proposition \ref{maibo} and the fact that for all $1 \leq i \leq n$ the density functions 
$f_i$ and $g_i$ belong to $\mathbb{L}_2(R)$. \hfill $\Box$\\

\noindent {\it{Proof of Proposition \ref{prop3}}}: 
\noindent Let us consider the variance of $T_j$. 
\noindent For all $(i_1,i_2)$, let $h_j(i_1,i_2)$ denote   the quantity
$$ h_j\left( i_1, i_2 \right) = \sum_k \sum_{l=1}^M \left( a_l(i_1)
\phi_{jk}(Y_{i_1}) - b_l(i_1) \phi_{jk}(Z_{i_1}) \right)
\left( a_l(i_2)
\phi_{jk}(Y_{i_2}) - b_l(i_2) \phi_{jk}(Z_{i_2}) \right).$$ 
\noindent The variance of $T_j$ satisfies
{\small{\begin{eqnarray*}
n^4 \ \mathbb{V}ar_{\overrightarrow{p},\overrightarrow{q}}(T_j)&=&\mathbb{V}ar_{\overrightarrow{p},\overrightarrow{q}}\left(\sum_{i_1\not= i_2}h_j(i_1,i_2)\right)\\
&=&\sum_{i_1\not= i_2, i_3\not= i_4}\mathbb{C}ov\left(h_j(i_1,i_2),h_j(i_3,i_4)\right)\\
&=&\sum_{i_1\not= i_2}\mathbb{V}ar\left(h_j(i_1,i_2)\right)+\sum_{i_1\not= i_2}\mathbb{C}ov\left(h_j(i_1,i_2),h_j(i_2,i_1)\right)\\
&&+\sum_{i_1\not= i_2 \not= i_3}\mathbb{C}ov\left(h_j(i_1,i_2),h_j(i_1,i_3)\right) +\sum_{i_1\not= i_2 \not= i_3}\mathbb{C}ov\left(h_j(i_1,i_2),h_j(i_2,i_3)\right)\\
&& + \sum_{i_1\not= i_2 \not= i_3}\mathbb{C}ov\left(h_j(i_1,i_2),h_j(i_3,i_1)\right)+\sum_{i_1\not= i_2 \not= i_3}\mathbb{C}ov\left(h_j(i_1,i_2),h_j(i_3,i_2)\right)\\
&&+\sum_{i_1\not= i_2 \not= i_3\not= i_4}\mathbb{C}ov\left(h_j(i_1,i_2),h_j(i_3,i_4)\right)\\
&=&\sum_{u=1}^7A_i.
\end{eqnarray*}}}

\noindent Using independence arguments, 
$$A_7=\sum_{i_1\not= i_2 \not= i_3\not= i_4}\mathbb{C}ov\left(h_j(i_1,i_2),h_j(i_3,i_4)\right)=0.$$

\noindent We are still required to bound for the quantities $A_i$ $(1 \leq i \leq 6).$ Since the ways to bound $A_1$ and $A_2$ (resp. $A_3,$ $A_4,$ $A_5$ and $A_6$) are similar, we will only bound $A_1$ and $A_3$. Such bounds are given in lemmas \ref{2ind} and \ref{3ind}. The proof of proposition \ref{prop3} is a direct consequence of lemmas \ref{2ind} and \ref{3ind} by taking $C_{_T}=2\ \bar{C}_{_T} \ \vee \ 4 \ \tilde{C}_{_T}$.\hfill
$\Box$ \\

\subsection{Proofs of Theorems}
\noindent{\it{Proof of Theorem \ref{theoup}.}}\\
\noindent Let us fix $0<\gamma<1$ and $s>0$. Under the null hypothesis, we use directly the well-known
Bienayme-Chebyshev inequality.
\begin{eqnarray*}
\mathbb{P}_{\overrightarrow{p},\overrightarrow{p}}\left( \Delta_s^* = 1 \right) & = & 
\mathbb{P}_{\overrightarrow{p},\overrightarrow{p}}\left( T_{j_n} > t_n \right) \nonumber \\
& \leq & \mathbb{P}_{\overrightarrow{p},\overrightarrow{p}}\left( T_{j_n}- \mathbb{E}(T_{j_n}) > t_n-\frac{8LMR^2}{Kn} \right) \nonumber \\
 & \leq & \frac {\mathbb{V}ar_{\overrightarrow{p},\overrightarrow{p}}(T_{j_n})} {\left(t_n-\frac{8LMR^2}{Kn}\right)^2} \nonumber \\
 & \leq & \frac{C_{_T} \ M^2 \ 2^{j_n}}{n^2 \ K^2 \left(t-\frac{8LMR^2}{K}\right)^2r_n^4}.
\label{nullNA}
\end{eqnarray*}
The last inequality is obtained using remark \ref{rem2}. According to the choices of the level $j_n$ and the threshold $t_n$, we have 
$$\frac{C_{_T} \ M^2 \ 2^{j_n}}{n^2 \ K^2\ \left(t-\frac{8LMR^2}{K}\right)^2r_n^4}
\leq \frac {2 C_{_T} \ M^2} {K^2\left(t-\frac{8LMR^2}{K}\right)^2}.$$ Then
$$\mathbb{P}_{\overrightarrow{p},\overrightarrow{p}}\left( \Delta_s^* = 1 \right)\leq \frac{\gamma}{2}.$$

\noindent Under the alternative, we use the expectation of the test statistic
and some approximation argument. The second type error is
\begin{eqnarray*}
\mathbb{P}_{\overrightarrow{p},\overrightarrow{q}}\left( \Delta_s^* = 0 \right) & = &
\mathbb{P}_{\overrightarrow{p},\overrightarrow{q}}\left( - T_{j_n} + \mathbb{E}_{\overrightarrow{p},\overrightarrow{q}}(T_{j_n}) \geq
- t_n + \mathbb{E}_{\overrightarrow{p},\overrightarrow{q}}(T_{j_n}) \right).
\end{eqnarray*}
The wavelet expansion in the Besov body $\mathcal{B}^s_{2, \infty}$ leads to
\begin{eqnarray*}
\mathbb{E}_{\overrightarrow{p},\overrightarrow{q}}(T_{j_n})  - t_n & = & 
\sum_{l=1}^M \left\| p_l - q_l \right\|_2^2
 - \sum_{l=1}^M \sum_{j \geq j_n} \sum_k \left( \int_{\mathbb{R}}
( p_l - q_l) \psi_{jk} \right)^2 \\
 & & \mathop{-} \frac{1}{n^2}
\sum_{l=1}^{M}\sum_k\sum_{i=1}^n\left(\int_{\mathbb{R}}\left(a_l(i)f_i-b_l(i)g_i\right)\phi_{j_nk}\right)^2 - t_n \\
 & \geq & \sum_{l=1}^M \left\| p_l - q_l \right\|_2^2
 - M \ R \ 2^{-2 j_n s} \mathop{-} \frac{8LMR^2}{Kn}- t_n. \\
  & \geq & \frac{1}{2}\sum_{l=1}^M \left\| p_l - q_l \right\|_2^2
 - M \ R \ 2^{-2 j_n s} - t_n,\\
 \end{eqnarray*} 
\noindent for any $n$ large enough.\\

\noindent  As a consequence, applying the Bienayme-Chebychev inequality leads
to
\begin{eqnarray*}
\lefteqn{\mathbb{P}_{\overrightarrow{p},\overrightarrow{q}}\left( - T_{j_n} + \mathbb{E}_{f,g}(T_{j_n}) \geq
- t_n + \mathbb{E}_{f,g}(T_{j_n}) \right)} \\
  & \leq & 
\frac{C_{_T} M^2\left(2^{j_n}+n\begin{displaystyle}\sum_l\end{displaystyle}\|p_l-q_l\|_2^2+ \sqrt{2^{j_n}n}\begin{displaystyle}\sum_l\end{displaystyle}\|p_l-q_l\|_2\right)}{n^2 \ K^2 \ \left(\begin{displaystyle}\frac{1}{2}\sum_{l=1}^M \end{displaystyle}\left\| p_l - q_l \right\|_2^2
 - M \ R \ 2^{-2 j_n s} - t_n\right)^2}.\end{eqnarray*}

\noindent The choice of $j_n$ and the fact that the functions are in the alternative entail
the following upper bound
\begin{eqnarray*}
\mathbb{P}_{\overrightarrow{p},\overrightarrow{q}}\left( \Delta_s^* = 0 \right) & \leq &
\frac{C_{_T}M^2\left(2^{j_n}+n\begin{displaystyle}\sum_l\end{displaystyle}\|p_l-q_l\|_2^2+ \sqrt{2^{j_n}n}\begin{displaystyle}\sum_l\end{displaystyle}\|p_l-q_l\|_2\right)}{K^2n^2\left(\frac{1}{2}\begin{displaystyle}\sum_{l=1}^M \end{displaystyle}\left\| p_l - q_l \right\|_2^2
 - M \ R \ 2^{-2 j_n s} - t \ r_n^2\right)^2}.\end{eqnarray*}

\noindent According to the choices of $j_n,$ and $r_n,$ one gets for  $n$ large enough:
\begin{eqnarray*}
\mathbb{P}_{\overrightarrow{p},\overrightarrow{q}}\left( \Delta_s^* = 0 \right) & \leq &
\frac{C_{_T}M^2\left(2^{j_n}+n\begin{displaystyle}\sum_l\end{displaystyle}\|p_l-q_l\|_2^2+ \sqrt{2^{j_n}n}\begin{displaystyle}\sum_l\end{displaystyle}\|p_l-q_l\|_2\right)}{n^2 \ K^2 \ \left(\frac{1}{2}
-\frac{R}{C^2}-\frac{t}{MC^2}\right)^2\left(\begin{displaystyle}\sum_{l=1}^M \end{displaystyle}\left\| p_l - q_l \right\|_2^2\right)^2}\\
 &\leq &\frac{3C_{_T}}{\left(\frac{1}{2}
-\frac{R}{C^2}-\frac{t}{MC^2}\right)^2K^2C^4}.
 \end{eqnarray*}
\noindent For all $C > C_\gamma,$ we finally obtain
 $$\mathbb{P}_{\overrightarrow{p},\overrightarrow{q}}\left( \Delta_s^* = 0 \right) \leq \frac{\gamma}{2}.$$

\noindent The results on the first-type and second-type errors show that if $C>C_\gamma$
the sum of the errors is less than $\gamma$. 
Therefore the upper bound is proved. \hfill $\Box$\\

\noindent {\it{Proof of Theorem \ref{theolow}.}} \\
Let $\gamma \in \left]0, \ 1 \right[,$ $C>0$ and $C_1>0$.
We define
\begin{eqnarray*}
\tilde{\Theta}_1\left( R,C,C_1,n,s\right) & =  & \Big\{ 
\left( \overrightarrow{p}, \overrightarrow{q} \right) :
\forall u \in \{1, \dots, M\}, p_u-q_u \in \mathcal{B}^s_{2,\infty}(R), \\ 
& & \quad \exists u \in \{1, \dots, M\}, \left(p_u,q_u\right) \in \tilde{\Lambda}_n(R,C,C_1) \Big\},
\end{eqnarray*} 
where $\tilde{\Lambda}_n(R,C,C_1)$ is defined in (\ref{spaceinf}).
It is well-known that 
{\small{\begin{eqnarray*}
&& \inf_{\Delta} \left( \sup_{(\overrightarrow{p},\overrightarrow{q}) \in \Theta_0(R)}\mathbb{P}_{\overrightarrow{p},\overrightarrow{q}}(\Delta=1) + \sup_{(\overrightarrow{p},\overrightarrow{q})  \in \Theta_1(R,C,n,s)} \mathbb{P}_{\overrightarrow{p},\overrightarrow{q}}(\Delta=0)\right)\\ &\geq&  \inf _{\Delta}
\left( \sup_{(\overrightarrow{p},\overrightarrow{q}) \in \Theta_0(R)}\mathbb{P}_{\overrightarrow{p},\overrightarrow{q}}(\Delta=1) + \sup_{(\overrightarrow{p},\overrightarrow{q})  \in \tilde{\Theta}_1(R,C,C_1,n,s)} \mathbb{P}_{\overrightarrow{p},\overrightarrow{q}}(\Delta=0)\right)\\
&\geq& 1 - \frac 1 2 \left\| \mathbb{P}_{\overrightarrow{p},\overrightarrow{p}} 
- \mathbb{P}_\pi \right\|,
\end{eqnarray*}}}
where $\left\| . \right\|$ is the $\mathbb L_1$-
distance and $\pi$ is an a priori probability measure on the set
$\Lambda_n(R,C)$.
\noindent First we define the probability measure $\pi$ and its support. Let  
$\theta=\left( \theta_1, \ldots, \theta_M \right)$ denote an eigenvector associated with the smallest 
eigenvalue of $\Sigma \Sigma^\star$ - which is $Kn$ according to
HYP-1 - such that
 $\left\| \theta \right\|_2 = 1$. \\
 
 \noindent Recall that here $j_n$ is the same as the one
 defined in theorem \ref{theoup}. Let $\mathcal{T}$ be the subset of $\mathbb{Z}$ containing every integer $k$ satisfying the following  properties
 
\begin{itemize} 
\item $k \in \mathcal{T} \Longrightarrow \left[\frac {k-L} {2^{j_n}},
\frac {k+L} {2^{j_n}} \right[ \subset [0,1[$;
\item $(k,k') \in \mathcal{T}\times \mathcal{T} \hbox{ with } k\not=k'\Longrightarrow \left[\frac {k-L} {2^{j_n}},
\frac {k+L} {2^{j_n}} \right[ \cap \left[\frac {k'-L} {2^{j_n}},
\frac {k'+L} {2^{j_n}} \right[ = \emptyset$.
\end{itemize}
The cardinal of $\mathcal{T}$ is clearly equal to $T=\lfloor \frac {2^{j_n-1}} L \rfloor$ 
and we denote its elements
$k_1, \ldots, k_T$.
The following
parametric family of functions is considered
$$q_{l,\zeta} (z) = p_l(z) +
2^{s+1}C  \sqrt{ML} 
\ \theta_l \ \sum_{k \in \mathcal{T}}  \zeta_k 
2^{-j_n s - \frac{j_n}{2}} \psi_{j_nk}(z),$$
where  $\zeta_k = +1$ or $-1$.

\noindent Remark that $\zeta_k$ does not depend on the index $l$.
Therefore the density of $Z_i$ is
$$ g_{i,\zeta}(z) = \sum_{l=1}^M \sigma_l(i)  \sqrt{ML} \ \theta_l \ 2^{s+1} C \
\sum_{k \in \mathcal{T}}  \zeta_k 2^{-j_n s - \frac{j_n}{2}} 
\psi_{j_nk}(z) + \sum_{l=1}^M \sigma_l(i)p_l(z).$$
The probability measure $\pi$ is such that the  $\zeta_k$'s are 
independent Rademacher random variables with
parameter $\frac 1 2$. 

\noindent The function $q_{l,\zeta}$ is a density. Indeed,
for $n$ large, $q_{l,\zeta}$ is non-negative.
Moreover, as $\psi_{j_nk}$ is a wavelet, we have
$\int \psi_{j_nk} = 0$ and therefore $\int q_{l,\zeta} = 1$.
If $C<\sqrt{R / M 2^{2s+2}}$, then $q_{l,\zeta} - p_l$ belongs to
the ball of the Besov body ${\mathcal B}^s_{2,\infty}(R)$.
There exists $l$ such that
$$ M \theta_l^2 \geq 1 \quad \mbox{ and } \quad
\left\| p_l - q_{l,\zeta} \right\|_2^2 =
T L M C^2 2^{2+2s-2j_ns-j_n} \theta_l^2 \geq 
C^2 \ n^{-\frac {4s} {4s+1}}.$$  

Therefore the probability measure $\pi$ is solely concentrated
on the alternative.

\noindent It is well-known that the ${\mathbb L}_1$ distance can
be bounded by the $\mathbb{L}_2$ distance. We have
\begin{eqnarray}
\left\| \mathbb{P}_{\overrightarrow{p},\overrightarrow{p}} - \mathbb{P}_\pi \right\| & \leq & 
\sqrt{\mathbb{E}_{\overrightarrow{p},\overrightarrow{p}}\left[
\left( \frac {d{\mathbb P}_\pi} {d{\mathbb P}_{\overrightarrow{p},\overrightarrow{p}}} \right)^2 \right] - 1}
 \nonumber \\
 & = &    \sqrt{\mathbb{E}_{\overrightarrow{p},\overrightarrow{p}}\left[
\left( \mathbb{E}_\pi\left(   
\prod_{i=1}^n \frac {g_{i,\zeta}(Z_i)}
{g_i(Z_i)} \right) \right)^2 \right] - 1}.\label{infi}
\end{eqnarray}

\noindent Therefore it suffices to evaluate the second-order moment of the likelihood ratio:
\begin{eqnarray*}
\lefteqn{\mathbb{E}_{\overrightarrow{p},\overrightarrow{p}}\left[
\left( \mathbb{E}_\pi\left(   
\prod_{i=1}^n \frac {g_{i,\zeta}(Z_i)}
{g_i(Z_i)} \right) \right)^2 \right]
} & & \\
 & & = \mathbb{E}_{\overrightarrow{p},\overrightarrow{p}}\left[ \left( \prod_{k \in \mathcal{T}} \int \prod_{i=1}^n
\left( 1 + 2^{s+1} C \sqrt{ML}  \ \zeta_k \ 2^{-j_ns-\frac{j_n}{2}} 
\frac {\psi_{j_nk}(Z_i)} {g_i(Z_i)} \sum_{l=1}^M  \theta_l  \sigma_l(i) 
 \right) \; 
d\pi(\zeta_1, \ldots, \zeta_T) \right)^2 \right]. 
\end{eqnarray*}
Let us introduce the following random variables
$$ \tilde{Z}_{ik} =  2^{s+1} C \sqrt{ML} \ 
 2^{-j_ns-\frac{j_n}{2}} \frac {\psi_{j_nk}(Z_i)} {g_i(Z_i)} \sum_{l=1}^M  \theta_l  \sigma_l(i).$$
\noindent We have
\begin{eqnarray*}
\lefteqn{\mathbb{E}_{\overrightarrow{p},\overrightarrow{p}}\left[ \left( \prod_{k \in \mathcal{T}} \int \prod_{i=1}^n
\left( 1 + 2^{s+1} C \sqrt{ML} \ \zeta_k \
2^{-j_ns-\frac{j_n}{2}} \frac {\psi_{j_nk}(Z_i)} {g_i(Z_i)} 
\sum_{l=1}^M  \theta_l  \sigma_l(i) \right) \; 
d\pi(\zeta_1, \ldots, \zeta_T) \right)^2 \right]} &  & \\
 & = & 
\mathbb{E}_{\overrightarrow{p},\overrightarrow{p}}\left[ \prod_{k \in \mathcal{T}} \frac 1 4 \left[ \prod_{i=1}^n
\left( 1 +  \tilde{Z}_{ik} \right) 
+ \prod_{i=1}^n \left( 1 -  \tilde{Z}_{ik} \right) \right]^2 \right]  \\
 & = & \mathbb{E}_{\overrightarrow{p},\overrightarrow{p}} \left[ \prod_{k \in \mathcal{T}} \frac 1 4 \left(
\prod_{i=1}^n (1+2 \tilde{Z}_{ik}+ \tilde{Z}_{ik}^2) + \prod_{i=1}^n (1-2 \tilde{Z}_{ik}+ \tilde{Z}_{ik}^2)
+ 2 \prod_{i=1}^n (1- \tilde{Z}_{ik}^2) \right) \right] \\
 & = & \mathbb{E}_{\overrightarrow{p},\overrightarrow{p}} \Big[ \prod_{k \in \mathcal{T}} \frac 1 4
\big\{ 2 \prod_{i=1}^n (1+ \tilde{Z}_{ik}^2) + 2 \prod_{i=1}^n
(1- \tilde{Z}_{ik}^2) \\
 & & \mathop{+} \sum_{i=1}^n \tilde{Z}_{ik} h_i(\tilde{Z}_{1k},\ldots, \tilde{Z}_{i-1,k},\tilde{Z}_{i+1,k}, \ldots,
\tilde{Z}_{nk}) \big\} \Big] \\
 & = & \mathbb{E}_{\overrightarrow{p},\overrightarrow{p}} \Big[ \prod_{k \in \mathcal{T}} \frac 1 2 \left(
\prod_{i=1}^n(1+ \tilde{Z}_{ik}^2) + \prod_{i=1}^n (1- \tilde{Z}_{ik}^2) \right)  \\
 & & \!\!\!\!\!\!\! \mathop{+} \sum_{r=1}^T \sum_{i=1}^n \tilde{Z}_{ik_r} \tilde{h}(\tilde{Z}_{1k_1}, \ldots, \tilde{Z}_{n,k_{r-1}},
\tilde{Z}_{1k_r}, \ldots,\tilde{Z}_{i-1,k_r}, \tilde{Z}_{i+1,k_r},\ldots, \tilde{Z}_{nk_r},\tilde{Z}_{1,k_{r+1}}, \ldots, \tilde{Z}_{nk_T} \Big],
\end{eqnarray*}
where the functions $h_i$ and $\tilde{h}_i$ are sums of
products of their arguments.
As $\mathbb{E}_{\overrightarrow{p},\overrightarrow{p}}( \tilde{Z}_{ik})=0$ and $\tilde{Z}_{ik} \tilde{Z}_{ik'} = 0$ for
$k \neq k'$, the last term vanishes. Thus we are only interested
in the first term.

\noindent Define for all $k \in \mathcal{T}$: 
\begin{eqnarray*}
h_l(k) & = & \sum_{1 \leq i_1 < i_2 < \ldots < i_l \leq n} \tilde{Z}_{i_1k}^2 \tilde{Z}_{i_2k}^2 \ldots
 \tilde{Z}_{i_lk}^2, \\
h_0(k) & = & 2.
\end{eqnarray*}
Then, we have
\begin{eqnarray*}
\mathbb{E}_{\overrightarrow{p},\overrightarrow{p}} \left[ \prod_{k \in \mathcal{T}} \frac 1 2 \left(
\prod_{i=1}^n(1+ \tilde{Z}_{ik}^2) + \prod_{i=1}^n (1- \tilde{Z}_{ik}^2) \right) \right] 
& = & \mathbb{E}_{\overrightarrow{p},\overrightarrow{p}} \left[ \left( \frac 1 2 \right)^T \prod_{k \in \mathcal{T}}
\left( \mathop{\sum_{l=0}}_{l \mbox{ even}}^n h_l(k) \right) \right] \\
 & = & \mathop{\sum_{l_1, \ldots, l_T=0}^n}_{l_1,\ldots,l_T
 \mbox{ even}} \left( \frac 1 2 \right)^T \mathbb{E}_{\overrightarrow{p},\overrightarrow{p}} \left(\prod_{r=1}^T
h_{l_r}(k_r) \right) 
\end{eqnarray*}
\begin{eqnarray*}
 & \leq & \mathop{\sum_{l_1, \ldots, l_T=0}^n \left( \frac 1 2 \right)^T}_{ l_1,\ldots, l_T \mbox{ even}} 
\prod_{r=1}^T \mathbb{E}_{\overrightarrow{p},\overrightarrow{p}}( h_{l_r}(k_r))  \\
 & \leq & \prod_{k \in \mathcal{T}} \frac 1 2 \left( \mathop{\sum_{l=0}}_{l \mbox{ even}}^n 
\mathbb{E}_{\overrightarrow{p},\overrightarrow{p}} \left[ h_l(k) \right] \right) \\
 & \leq & \prod_{k \in \mathcal{T}} \frac 1 2 \left( \mathop{\sum_{l=0}}_{l \mbox{ even}}^n 
\sum_{1 \leq i_1 < \ldots < i_l \leq n}^n \mathbb{E}_{\overrightarrow{p},\overrightarrow{p}}\left[ \tilde{Z}_{i_1k}^2 \right]
\ldots \mathbb{E}\left[ \tilde{Z}_{i_lk}^2 \right] \right) \\
 & \leq & \prod_{k \in \mathcal{T}} \frac 1 2 \left( \prod_{i=1}^n 
\left( 1 + \mathbb{E}_{\overrightarrow{p},\overrightarrow{p}}\left[ \tilde{Z}_{ik}^2 \right] \right) + \prod_{i=1}^n
\left( 1 - \mathbb{E}_{\overrightarrow{p},\overrightarrow{p}}\left[ \tilde{Z}_{ik}^2 \right] \right) \right) \\
 & \leq & \prod_{k \in \mathcal{T}} \cosh \left( \sum_{i=1}^n \mathbb{E}_{\overrightarrow{p},\overrightarrow{p}}\left( 
\tilde{Z}_{ik}^2 \right) \right) \\
 & \leq & \exp\left( \frac{1}{2} \sum_{k \in \mathcal{T}} \left( \sum_{i=1}^n
\mathbb{E}_{\overrightarrow{p},\overrightarrow{p}}\left( \tilde{Z}_{ik}^2 \right) \right)^2 \right). 
\end{eqnarray*}

\noindent Each $ \mathbb{E}_{\overrightarrow{p},\overrightarrow{p}}\left( \tilde{Z}_{ik}^2 \right)$ is bounded as follows,
$$ \mathbb{E}_{\overrightarrow{p},\overrightarrow{p}}\left( \tilde{Z}_{ik}^2 \right) \leq
2^{2s+2-2j_ns-j_n} \frac {C^2} {C_1} 
M L \left( \sum_{l=1}^M  \theta_l  \sigma_l(i) \right)^2.$$
Therefore this bound entails
{\small{\begin{eqnarray}
\exp\left( \frac{1}{2} \sum_{k \in \mathcal T} \left( \sum_{i=1}^n
\mathbb{E}_{\overrightarrow{p},\overrightarrow{p}}\left( \tilde{Z}_{ik}^2 \right) \right)^2 \right) & \leq & 
\exp\left( \frac{1}{2} \sum_{k \in \mathcal T} 2^{4s+4}  C^4 2^{-4j_ns-2j_n} \frac {L^2 M^2} {C_1^2} \left( \sum_{i=1}^n
\sum_{l,m=1}^M \theta_l \theta_m \sigma_l(i) \sigma_m(i) \right)
^2   \right) \nonumber \\
&\leq&\exp\left( \frac{1}{2} \sum_{k \in \mathcal T} 2^{4s+4}  C^4 2^{-4j_ns-2j_n} \frac {L^2 M^2} {C_1^2} \left(\theta^\star  \Gamma_n' \theta
\right)^2\right) \nonumber \\
&=&\exp\left( \sum_{k \in \mathcal T} 2^{4s+3}  C^4 2^{-4j_ns-2j_n} \frac {L^2 M^2} {C_1^2} \left(Kn\right)^2\right) \nonumber \\
 & \leq & \exp\left( 2^{4s+2} M^2K^2  \ \frac {LC^4 }
{C_1^2} \right). \label{finallowerbound}
\end{eqnarray}}}

\noindent Inequalities (\ref{infi}) and (\ref{finallowerbound}) lead to
$$\left\| \mathbb{P}_{\overrightarrow{p},\overrightarrow{p}} - \mathbb{P}_\pi \right\| 
\leq \sqrt{\exp\left( 2^{4s+2} M^2K^2  \ \frac {LC^4 }
{C_1^2} \right)-1}.$$
The choice of any constant $C$ such that $C< c_\gamma$ entails that
 the left-hand side of (\ref{finallowerbound}) is strictly smaller than $
2(1-\gamma)$.

\hfill $\Box$\\

\section{Appendix}
\noindent This section contains the technical lemmas used in the proofs
of the main results.\\

\begin{lemma}\label{lem0}
\begin{eqnarray}
\sum_{l=1}^M\sum_{i=1}^na_l^2(i)&\leq &  \frac{Mn}{K},\label{tr1}\\
\sum_{l=1}^M\sum_{i=1}^nb_l^2(i)&\leq &  \frac{Mn}{K}.\label{tr2}
\end{eqnarray}

\end{lemma}

\noindent {\it{Proof of Lemma \ref{lem0}}}:  \\
The proofs of (\ref{tr1}) and (\ref{tr2}) are identical, that's why we only prove (\ref{tr1}).  \noindent Let $\lambda_{min}(\Gamma_n)$ be the smallest non negative eigenvalue of the matrix $\Gamma_n$. Let 
$A=(A)_{1 \leq j \leq n, 1 \leq l \leq M}$ denote the $(n \times M)$ matrix  with coefficients $A_{j,l}=a_{l}(j)$. Since the matrix $AA^*$ has at most $M$ non negative eigenvalues, we have
 \begin{eqnarray}\sum_{l=1}^M \sum_{i=1}^n a_l^2(i)=trace(AA^*)\leq M \ \lambda_{max}(AA^*).\label{sed1}\end{eqnarray}
Clearly, the following implication holds
 $$\lambda \hbox{ is a non negative eigenvalue of } AA^* \Longrightarrow n^2\lambda^{-1} \hbox{ is an eigenvalue of } \Gamma_n.$$
\noindent So
\begin{eqnarray}\lambda_{max}(AA^*)\leq \frac{n^2}{\lambda_{min}(\Gamma_n)}.\label{sed2}
\end{eqnarray}
Lemma \ref{lem0} is proved by inequalities (\ref{sed1}) and (\ref{sed2}) and 
under HYP-1. 
 \hfill
$\Box$ \\

\begin{lemma}\label{lem1}
For all $(j,k) \in \mathbb{Z}\times \mathbb{Z}$, let us put $$I_{jk}=
\left[\frac{k-L}{2^j},\frac{k+L}{2^j}\right[.$$
Then for any fixed $(j,k)$
$$Card\{k' \in \mathbb{Z}: \ I_{jk}\cap I_{jk'} \not= \emptyset\}\leq 4L.$$
\end{lemma}

\noindent {\it{Proof of Lemma \ref{lem1}}}:  \\
Clearly,
$I_{jk}\cap I_{jk'} = \emptyset \iff k'-L\geq k+L\ \hbox{ or } \ k'+L\leq k-L.$ \\ 
Hence, $I_{jk}\cap I_{jk'} \not= \emptyset \iff k-2L<k'<k+2L.$ \\
As a consequence, we have
 $$Card\{k' \in \mathbb{Z}: \ I_{jk}\cap I_{jk'} \not= \emptyset\}\leq 4L.$$
 \hfill
$\Box$ \\

\begin{lemma}\label{lem2}
For any function $h \in L_1(\mathbb{R})$ 
\begin{eqnarray*}\sum_k \int_{I_{jk}}|h(x)|dx \leq 2L \|h\|_1.
\end{eqnarray*}
\end{lemma}

\noindent {\it{Proof of Lemma \ref{lem2}}}: Let us define for any $h \in L_1(\mathbb{R})$ :
 $$p_{jk}(h)=\int_{I_{jk}}|h(x)|dx, \quad \forall j\in \mathbb{N},\ \ \forall k \in \mathbb{Z}.$$

\noindent Judging from the definition of the intervals $I_{jk}$, we easily prove that for any $j\in \mathbb{N}$,
\begin{eqnarray*}
\sum_k p_{jk}(h)=\sum_{u=1}^{2L}\sum_{i \in \mathbb{Z}} p_{j,2Li+u}(h) \leq
\sum_{u=1}^{2L}\int_\mathbb{R} |h(x)|dx=2L \|h\|_1.
\end{eqnarray*} \hfill
$\Box$ \\

\begin{lemma}\label{teclem1}
\noindent Let $W$ be either $Y$ or $Z$. For any $1 \leq i \leq n$  and any $(j,k)$,  we have
\begin{eqnarray*}
\left| \mathbb E\left( \phi_{jk}( W_{i} ) \right) \right| & \leq & 
\left( 2 L \ \sup_{l} (\left\| p_l \right\|_\infty \vee\left\| q_l \right\|_\infty) \right)^{\frac 1 2} 
2^{- \frac j 2}.
\end{eqnarray*}
\end{lemma}

\noindent {\it{Proof of Lemma \ref{teclem1}}}:

\noindent Using the Cauchy-Schwarz inequality, we obtain
\begin{eqnarray*}
\left| \mathbb E\left( \phi_{jk}( W_{i} ) \right) \right| & \leq& 
\left|\int\phi_{jk}\ f_i\right| \vee \left|\int\phi_{jk}\ g_i\right|\\
& \leq & \int|\phi_{jk}|\ \sup_{l} \left\| p_l \right\|_\infty \vee \int|\phi_{jk}|\ \sup_{l} \left\| q_l \right\|_\infty\\
&\leq&\left( 2 L \ \sup_{l} (\left\| p_l \right\|_\infty \vee\left\| q_l \right\|_\infty) \right)^{\frac 1 2} 
2^{- \frac j 2}.
\end{eqnarray*}
\hfill
$\Box$ \\

\begin{lemma}\label{teclem2}
\noindent Let $W$ be either $Y$ or $Z$ and $c$ be either $a$ or $b$. For any $1 \leq i \leq n$  and any $(j,k)$, the following inequalities hold
\begin{eqnarray*}
\sum_{k'}\left| \mathbb E\left( \phi_{jk}(W_{i}) \phi_{jk'}(W_{i}) \right)
\right| & \leq & 4L \ \sup_{l} (\left\| p_l \right\|_\infty \vee\left\| q_l \right\|_\infty), \\
\sup_l\left|\sum_{k}  \int \phi_{jk}  (p_l - q_l) \right| & \leq & 
4L \left\| \phi \right\|_\infty 2^{\frac j 2}, \\
\sup_l\left| c_l(i) \right| & \leq & \sqrt{n \sum_{l} \left< c_l, c_l \right>_n}.
\end{eqnarray*}
\end{lemma}

\noindent {\it{Proof of Lemma \ref{teclem2}}}:\\
\noindent Since the wavelets are compactly supported, for any fixed $k$ the sum over $k'$ has at most $4L$ terms which are non zeros (see lemma \ref{lem1}). So, the Cauchy-Schwarz inequality entails that

\begin{eqnarray*}
\sum_{k'}\left| \mathbb E\left( \phi_{jk}(W_{i}) \phi_{jk'}(W_{i}) \right)
\right| & \leq& \sum_{k'}\int \left|\phi_{jk}\right|\left|\phi_{jk'}\right| f_i \vee \sum_{k'}\int \left|\phi_{jk}\right|\left|\phi_{jk'}\right| g_i\\
& \leq& \sum_{k'}\left(\|f_i\|_\infty\int \left|\phi_{jk}\right|\left|\phi_{jk'}\right|\right)  \vee \sum_{k'}\left(\|g_i\|_\infty\int \left|\phi_{jk}\right|\left|\phi_{jk'}\right| \right)\\
& \leq& \left(\sup_l\left\| p_l \right\|_\infty \sum_{k'}\int \left|\phi_{jk}\right|\left|\phi_{jk'}\right|\right)  \vee \left(\sup_l\left\| q_l \right\|_\infty\sum_{k'}\int \left|\phi_{jk}\right|\left|\phi_{jk'}\right| \right)\\
& \leq & 4L \ \sup_{l} (\left\| p_l \right\|_\infty \vee\left\| q_l \right\|_\infty).
\end{eqnarray*}
We also have
\begin{eqnarray*}
\sup_l\left|\sum_{k}  \int \phi_{jk}  (p_l - q_l) \right| & \leq & 
 2^{\frac{j}{2}}\|\phi\|_\infty \sup_l \sum_{k}  \int_{I_{jk}}  |p_l - q_l| \\
 & \leq &  2L \left(\int p_l + \int q_l \right)\|\phi\|_\infty \ 2^{\frac{j}{2}}\\
 & = &  4L \|\phi\|_\infty \ 2^{\frac{j}{2}}.
\end{eqnarray*}

\noindent Clearly, for any $1 \leq i \leq n$,
\begin{eqnarray*}
\sup_l\left| c_l(i) \right| & \leq & \sup_l \sqrt{\sum_i c_l^2(i)}\\
& \leq & \sqrt{n \sum_{l} \left< c_l, c_l \right>_n}.
\end{eqnarray*}\hfill
$\Box$ \\

\begin{lemma}\label{teclem3}
\noindent Let $p_l,$ $q_l$ $p_{l'}$ and $q_{l'}$ be four probability 
densities in $\mathbb{L}_2$. Then, for any $j \in \mathbb{N}$

\begin{eqnarray*}
\sum_k\left(\int\phi_{jk}p_l-\int\phi_{jk}q_l\right)^2 &\leq& 2L\|p_l-q_l\|_2^2;\\
\sum_{k}\sum_{k'  : I_{jk} \cap I_{jk'} \neq \emptyset}
 \left| (\int \phi_{jk} p_l - \int \phi_{jk} q_l)(
\int \phi_{jk'} p_{l'} - \int \phi_{jk'} q_{l'}) \right| &\leq& 
4L^2 \left(\|p_l-q_l\|_2^2+\|p_{l'}-q_{l'}\|_2^2\right).
\end{eqnarray*}
\end{lemma}

\noindent {\it{Proof of Lemma \ref{teclem3}}}:\\
Using the Cauchy-Schwarz inequality, we have
\begin{eqnarray*}
\sum_k\left(\int\phi_{jk}p_l-\int\phi_{jk}q_l\right)^2 & \leq& \sum_k\int_{I_{jk}}(p_l-q_l)^2\\
&\leq& 2L\|p_l-q_l\|_2^2.
\end{eqnarray*}

\noindent Lemma \ref{lem2} entails that

\begin{eqnarray*}
&&\sum_{k}\sum_{k'  : I_{jk} \cap I_{jk'} \neq \emptyset}
 \left| (\int \phi_{jk} p_l - \int \phi_{jk} q_l)(
\int \phi_{jk'} p_{l'} - \int \phi_{jk'} q_{l'}) \right|\\
& \leq &
\frac 1 2 \left[\sum_{k}\sum_{k'  : I_{jk} \cap I_{jk'} \neq \emptyset}
 \left(\int \phi_{jk} p_l - \int \phi_{jk} q_l \right)^2 + \sum_{k}\sum_{k'  : I_{jk} \cap I_{jk'} \neq \emptyset}
 \left(\int \phi_{jk} p_{l'} - \int \phi_{jk} q_{l'} \right)^2\right]  \\  
&\leq& \frac 1 2 \left[4L\sum_k\int\phi_{jk}^2\int_{I_{jk}}(p_l-q_l)^2+
4L\sum_k\int\phi_{jk}^2\int_{I_{jk}}(p_{l'}-q_{l'})^2\right]\\
&\leq & \frac 1 2 \left(8L^2\|p_l-q_l\|_2^2+8L^2\|p_{l'}-q_{l'}\|_2^2\right)\\
&\leq& 4L^2 \left(\|p_l-q_l\|_2^2+\|p_{l'}-q_{l'}\|_2^2\right).
\end{eqnarray*}
\hfill
$\Box$ \\

\begin{lemma}\label{2ind}
There exists a constant $\bar{C}_{_T}=\bar{C}_{_T}(R,L,\|\phi\|_\infty)>0$ such that $$A_1:=\sum_{i_1\not= i_2}\mathbb{V}ar_{_{\overrightarrow{p},\overrightarrow{q}}}(h_j(i_1,i_2)) \leq \bar{C}_{_T} \frac{M^2}{K^2} \ 2^j \ n^2.$$
\end{lemma}

\noindent {\it{Proof of Lemma \ref{2ind}}}:\\
Let us evaluate each variance
\begin{eqnarray*} \mathbb{V}ar_{_{\overrightarrow{p},\overrightarrow{q}}}\left( h_j\left( i_1, i_2 \right) \right) &=& \mathbb{C}ov
\left( h_j(i_1,i_2), h_j(i_1,i_2) \right).
\end{eqnarray*}
We expand the covariance
\begin{eqnarray*}
\lefteqn{\mathbb{C}ov\big( 
\left( a_l(i_1) \phi_{jk}(Y_{i_1}) - b_l(i_1) \phi_{jk}(Z_{i_1}) \right)
\left( a_l(i_2) \phi_{jk}(Y_{i_2}) - b_l(i_2) \phi_{jk}(Z_{i_2}) \right), } \\
\lefteqn{  \left( a_{l'}(i_1) \phi_{jk'}(Y_{i_1}) - b_{l'}(i_1) \phi_{jk'}(Z_{i_1}) \right)
\left( a_{l'}(i_2) \phi_{jk'}(Y_{i_2}) - b_{l'}(i_2) \phi_{jk'}(Z_{i_2}) \right) \big)} 
 \\
 & = & 
\mathbb{C}ov\left(  a_l(i_1) \phi_{jk}(Y_{i_1}) 
 a_l(i_2) \phi_{jk}(Y_{i_2}),
 a_{l'}(i_1) \phi_{jk'}(Y_{i_1}) 
 a_{l'}(i_2) \phi_{jk'}(Y_{i_2}) \right) \\
 & & - \mathbb{C}ov\left(  a_l(i_1) \phi_{jk}(Y_{i_1}) 
 a_l(i_2) \phi_{jk}(Y_{i_2}), 
 a_{l'}(i_1) \phi_{jk'}(Y_{i_1}) b_{l'}(i_2) \phi_{jk'}(Z_{i_2}) \right) \\
 & & -  \mathbb{C}ov\left(  a_l(i_1) \phi_{jk}(Y_{i_1}) 
 a_l(i_2) \phi_{jk}(Y_{i_2}),
 b_{l'}(i_1) \phi_{jk'}(Z_{i_1})
 a_{l'}(i_2) \phi_{jk'}(Y_{i_2}) \right) \\
 & &  + \mathbb{C}ov\left(  a_l(i_1) \phi_{jk}(Y_{i_1}) 
 a_l(i_2) \phi_{jk}(Y_{i_2}),
 b_{l'}(i_1) \phi_{jk'}(Z_{i_1})
 b_{l'}(i_2) \phi_{jk'}(Z_{i_2}) \right)  \\
 & & - \mathbb{C}ov\left( a_l(i_1) \phi_{jk}(Y_{i_1})
 b_l(i_2) \phi_{jk}(Z_{i_2}),
 a_{l'}(i_1) \phi_{jk'}(Y_{i_1}) 
 a_{l'}(i_2) \phi_{jk'}(Y_{i_2}) \right) \\
 & & + \mathbb{C}ov\left( a_l(i_1) \phi_{jk}(Y_{i_1})
 b_l(i_2) \phi_{jk}(Z_{i_2}),
 a_{l'}(i_1) \phi_{jk'}(Y_{i_1}) 
 b_{l'}(i_2) \phi_{jk'}(Z_{i_2}) \right) \\
 & & + \mathbb{C}ov\left( a_l(i_1) \phi_{jk}(Y_{i_1})
 b_l(i_2) \phi_{jk}(Z_{i_2}),
 b_{l'}(i_1) \phi_{jk'}(Z_{i_1}) 
 a_{l'}(i_2) \phi_{jk'}(Y_{i_2}) \right) \\
 & & - \mathbb{C}ov\left( a_l(i_1) \phi_{jk}(Y_{i_1})
 b_l(i_2) \phi_{jk}(Z_{i_2}),
 b_{l'}(i_1) \phi_{jk'}(Z_{i_1}) 
 b_{l'}(i_2) \phi_{jk'}(Z_{i_2}) \right) \\
 & & -\mathbb{C}ov\left( b_l(i_1) \phi_{jk}(Z_{i_1})
 a_l(i_2) \phi_{jk}(Y_{i_2}),
 a_{l'}(i_1) \phi_{jk'}(Y_{i_1})
 a_{l'}(i_2) \phi_{jk'}(Y_{i_2}) \right) \\
 & & + \mathbb{C}ov\left( b_l(i_1) \phi_{jk}(Z_{i_1})
 a_l(i_2) \phi_{jk}(Y_{i_2}),
 a_{l'}(i_1) \phi_{jk'}(Y_{i_1})
 b_{l'}(i_2) \phi_{jk'}(Z_{i_2}) \right) \\
 & & + \mathbb{C}ov\left( b_l(i_1) \phi_{jk}(Z_{i_1})
 a_l(i_2) \phi_{jk}(Y_{i_2}),
 b_{l'}(i_1) \phi_{jk'}(Z_{i_1})
 a_{l'}(i_2) \phi_{jk'}(Y_{i_2}) \right) \\
 & & -  \mathbb{C}ov\left( b_l(i_1) \phi_{jk}(Z_{i_1})
 a_l(i_2) \phi_{jk}(Y_{i_2}),
 b_{l'}(i_1) \phi_{jk'}(Z_{i_1})
 b_{l'}(i_2) \phi_{jk'}(Z_{i_2}) \right) \\
 & & + \mathbb{C}ov\left( b_l(i_1) \phi_{jk}(Z_{i_1})
 b_l(i_2) \phi_{jk}(Z_{i_2}),
 a_{l'}(i_1) \phi_{jk'}(Y_{i_1})
 a_{l'}(i_2) \phi_{jk'}(Y_{i_2}) \right) \\
 & & - \mathbb{C}ov\left( b_l(i_1) \phi_{jk}(Z_{i_1})
 b_l(i_2) \phi_{jk}(Z_{i_2}),
 a_{l'}(i_1) \phi_{jk'}(Y_{i_1})
 b_{l'}(i_2) \phi_{jk'}(Z_{i_2}) \right) \\
 & & - \mathbb{C}ov\left( b_l(i_1) \phi_{jk}(Z_{i_1})
 b_l(i_2) \phi_{jk}(Z_{i_2}),
 b_{l'}(i_1) \phi_{jk'}(Z_{i_1})
 a_{l'}(i_2) \phi_{jk'}(Y_{i_2}) \right) \\
 & & + \mathbb{C}ov\left( b_l(i_1) \phi_{jk}(Z_{i_1})
 b_l(i_2) \phi_{jk}(Z_{i_2}),
 b_{l'}(i_1) \phi_{jk'}(Z_{i_1})
 b_{l'}(i_2) \phi_{jk'}(Z_{i_2}) \right).
\end{eqnarray*}
According to independence arguments, the following terms are clearly equal to zero:
$$
\mathbb{C}ov\left(  a_l(i_1) \phi_{jk}(Y_{i_1}) 
 a_l(i_2) \phi_{jk}(Y_{i_2}),
 b_{l'}(i_1) \phi_{jk'}(Z_{i_1})
 b_{l'}(i_2) \phi_{jk'}(Z_{i_2}) \right),$$
$$ \mathbb{C}ov\left( a_l(i_1) \phi_{jk}(Y_{i_1})
 b_l(i_2) \phi_{jk}(Z_{i_2}),
 b_{l'}(i_1) \phi_{jk'}(Z_{i_1}) 
 a_{l'}(i_2) \phi_{jk'}(Y_{i_2}) \right),$$
$$ 
 \mathbb{C}ov\left( b_l(i_1) \phi_{jk}(Z_{i_1})
 a_l(i_2) \phi_{jk}(Y_{i_2}),
 a_{l'}(i_1) \phi_{jk'}(Y_{i_1})
 b_{l'}(i_2) \phi_{jk'}(Z_{i_2}) \right),$$
$$ \mathbb{C}ov\left( b_l(i_1) \phi_{jk}(Z_{i_1})
 b_l(i_2) \phi_{jk}(Z_{i_2}),
 a_{l'}(i_1) \phi_{jk'}(Y_{i_1})
 a_{l'}(i_2) \phi_{jk'}(Y_{i_2}) \right).$$ 
The remaining terms can be split into two types: those involving
two different random variables and those involving 
three different random variables. Let us handle these two cases separately.
First, we consider the case with two different random variables.
 We need to bound terms such as
\begin{eqnarray*}
\lefteqn{\sum_{i_1 \not= i_2}\sum_{k,k'} \mathbb{C}ov\left(  a_l(i_1) \phi_{jk}(Y_{i_1}) 
 a_l(i_2) \phi_{jk}(Y_{i_2}),
 a_{l'}(i_1) \phi_{jk'}(Y_{i_1}) 
 a_{l'}(i_2) \phi_{jk'}(Y_{i_2}) \right)} \\
 & = &  \sum_{i_1 \not= i_2}\sum_{k,k'} a_l(i_1) a_l(i_2) a_{l'}(i_1) a_{l'}(i_2)
\mathbb E\left( \phi_{jk}(Y_{i_1}) \phi_{jk'}(Y_{i_1}) \right)
\mathbb E\left( \phi_{jk}(Y_{i_2}) \phi_{jk'}(Y_{i_2}) \right)
 \\
 & & - \sum_{i_1 \not= i_2}\sum_{k,k'} a_l(i_1) a_l(i_2) a_{l'}(i_1) a_{l'}(i_2)
\mathbb E\left( \phi_{jk}(Y_{i_1}) \right)
\mathbb E\left( \phi_{jk'}(Y_{i_1}) \right)
\mathbb E\left( \phi_{jk}(Y_{i_2}) \right)
\mathbb E\left( \phi_{jk'}(Y_{i_2}) \right).
\end{eqnarray*}
As the wavelets are compactly supported, we get  for any $(i_1,i_2)$,
\begin{eqnarray*}
&&\left| \sum_{k,k'} a_l(i_1) a_l(i_2) a_{l'}(i_1) a_{l'}(i_2)
\mathbb E\left( \phi_{jk}(Y_{i_1}) \phi_{jk'}(Y_{i_1}) \right)
\mathbb E\left( \phi_{jk}(Y_{i_2}) \phi_{jk'}(Y_{i_2}) \right) \right| \\
& \leq &  \|f_{i_1}\|_\infty \sum_{k,k'} |a_l(i_1) a_l(i_2) a_{l'}(i_1) a_{l'}(i_2)|\int|\phi_{jk}\phi_{jk'}|f_{i_2}\\
 & \leq &  2^{j+3} L^2 \left\| \phi \right\|_\infty^2 \sup_{l} (\left\| p_l \right\|_\infty \vee\left\| q_l \right\|_\infty)
|a_l(i_1) a_l(i_2) a_{l'}(i_1) a_{l'}(i_2)|.
\end{eqnarray*}
The second sum is much simpler to bound. According to lemma \ref{teclem1}
it can be bounded as follows
\begin{eqnarray*}
 \lefteqn{\left| \sum_{k,k'} a_l(i_1) a_l(i_2) a_{l'}(i_1) a_{l'}(i_2)
\mathbb E\left( \phi_{jk}(Y_{i_1}) \right)
\mathbb E\left( \phi_{jk'}(Y_{i_1}) \right)
\mathbb E\left( \phi_{jk}(Y_{i_2}) \right)
\mathbb E\left( \phi_{jk'}(Y_{i_2}) \right) \right|} \\
 & \leq & \sum_{k,k'} |a_l(i_1) a_l(i_2) a_{l'}(i_1) a_{l'}(i_2)|
\mathbb E\left( |\phi_{jk}(Y_{i_1})| \right)
\mathbb E\left( |\phi_{jk'}(Y_{i_1})| \right)
\left(\sqrt{2L} \  2^{-\frac{j}{2}} \right)^2 \ \sup_{l} (\left\| p_l \right\|_\infty \vee\left\| q_l \right\|_\infty) \\
 & =  &  L \ 2^{1-j}\sum_{k,k'} |a_l(i_1) a_l(i_2) a_{l'}(i_1) a_{l'}(i_2)|
\int_{I_{jk}} |\phi_{jk}| f_{{i_1}} 
\int_{I_{jk'}} |\phi_{jk'}| f_{{i_1}}
\left( \sup_{l} (\left\| p_l \right\|_\infty \vee\left\| q_l \right\|_\infty) \right) \\
 & \leq &  8 L^3 \left| a_l(i_1) a_l(i_2) a_{l'}(i_1) a_{l'}(i_2) \right| 
\left\| \phi \right\|_\infty^2
\sup_{l} (\left\| p_l \right\|_\infty \vee\left\| q_l\right\|_\infty).
\end{eqnarray*}
Let us now focus on the sums over $i_1, i_2, l$ and $l'$.
\begin{eqnarray*}
\sum_{i_1 \neq i_2}\sum_{l,\ l'} \left| a_l(i_1) a_{l}(i_2) a_{l'}(i_1) a_{l'}(i_2) \right| & \leq & 
\sum_{i_1, i_2} \sum_{l,\ l'} \frac 1 2 \left( a_l(i_1)^2 a_{l'}(i_2)^2 + a_{l'}(i_1)^2 a_{l}(i_2)^2 \right) \\
 & \leq & n^2 \sum_{l,\ l'} \left< a_l, a_l \right>_n
\ \left< a_{l'}, a_{l'}\right>_n \\\
 & \leq & \frac{M^2n^2}{K^2}.
\end{eqnarray*}
We see that this term behaves like $n^2$.
The three other terms featuring only two different random variables are handled in the same way.

\noindent Therefore it remains to evaluate the eight terms with three different random variables. For example, let 
us consider
$$ \mathbb{C}ov\left( a_l(i_1) \phi_{jk}(Y_{i_1}) a_l(i_2) \phi_{jk}(Y_{i_2}),
a_{l' }(i_1) \phi_{jk'}(Y_{i_1}) a_{l'}(i_2) \phi_{jk'}(Z_{i_2}) \right),$$
and let us omit for a moment
the sums over $i_1, i_2, k, k', l$ and $l'$. 
The covariance can be expanded as
\begin{eqnarray*}
 \mathbb{C}ov\left(  \phi_{jk}(Y_{i_1})  \phi_{jk}(Y_{i_2}),
 \phi_{jk'}(Y_{i_1})  \phi_{jk'}(Z_{i_2}) \right) 
 & = & 
 \mathbb E\left( \phi_{jk'}(Z_{i_2}) \right)
\mathbb E\left( \phi_{jk}(Y_{i_2}) \right)\mathbb{C}ov\left(\phi_{jk}(Y_{i_1}),\phi_{jk'}(Y_{i_1})\right).
\end{eqnarray*}
When we add the sums over $k$ and $k'$, the second term is exactly handled as the second term above
in the case of two different random variables. Thus, it remains to consider the first summand. As above, the
compactness of the wavelet entails that
\begin{eqnarray*}
&&\left| \sum_{k,k'} \mathbb E\left( \phi_{jk'}(Z_{i_2}) \right)
\mathbb E\left( \phi_{jk}(Y_{i_2}) \right)\mathbb{C}ov\left(\phi_{jk}(Y_{i_1}),\phi_{jk'}(Y_{i_1})\right) \right|  \\
& \leq & \sum_{k,k'}|\mathbb E\left( \phi_{jk'}(Z_{i_2}) \right)
\mathbb E\left( \phi_{jk}(Y_{i_2}) \right)\mathbb E\left( \phi_{jk}(Y_{i_1}) \right)
\mathbb E\left( \phi_{jk'}(Y_{i_1}) \right)|\\
&& \quad \quad +\sum_{k,k'}|\mathbb E\left( \phi_{jk'}(Z_{i_2}) \right)
\mathbb E\left( \phi_{jk}(Y_{i_2}) \right)\mathbb E\left( \phi_{jk}(Y_{i_1})\phi_{jk'}(Y_{i_1}) \right)|\\
&=& A_{11}+A_{12},
\end{eqnarray*}

\noindent According to lemmas \ref{teclem1} and \ref{teclem2}, we have
\begin{eqnarray*}
A_{11} &=&\sum_{k,k'}|\mathbb E\left( \phi_{jk'}(Z_{i_2}) \right)
\mathbb E\left( \phi_{jk}(Y_{i_2}) \right)\mathbb E\left( \phi_{jk}(Y_{i_1}) \right)
\mathbb E\left( \phi_{jk'}(Y_{i_1}) \right)|\\
&\leq&\left(2^{1-j}L\sup_{l} (\left\| p_l \right\|_\infty \vee\left\| q_l \right\|_\infty)\right)\sum_{k}\int|\phi_{jk}|g_{i_2}\sum_{k'}\int|\phi_{jk'}|f_{i_2}\\
&\leq& 8L^3 \left\| \phi \right\|_\infty^2\sup_{l} (\left\| p_l \right\|_\infty \vee\left\| q_l \right\|_\infty)
\end{eqnarray*}
and
\begin{eqnarray*}
A_{12} &=&\sum_{k,k'}|\mathbb E\left( \phi_{jk'}(Z_{i_2}) \right)
\mathbb E\left( \phi_{jk}(Y_{i_2}) \right)\mathbb E\left( \phi_{jk}(Y_{i_1})\phi_{jk'}(Y_{i_1}) \right)|\\
&\leq & \left(2^{1-j}L\sup_{l} (\left\| p_l \right\|_\infty \vee\left\| q_l \right\|_\infty)\right)\sum_{k,k'}|\mathbb E\left( \phi_{jk}(Y_{i_1})\phi_{jk'}(Y_{i_1}) \right)|\\
&\leq & 4L \left(2^{1-j}L\sup_{l} (\left\| p_l \right\|_\infty \vee\left\| q_l \right\|_\infty)\right)2^{\frac{j}{2}}\|\phi\|_\infty\sum_{k}\int|\phi_{jk}|f_{i_1}\\
&\leq & 8L^2\left(2L\sup_{l} (\left\| p_l \right\|_\infty \vee\left\| q_l \right\|_\infty)\right)\|\phi\|_\infty^2\int f_{i_1}\\
& \leq & 16 L^3\left\|\phi\right\|_\infty^2\sup_{l} (\left\| p_l \right\|_\infty \vee\left\| q_l \right\|_\infty).
\end{eqnarray*}

\noindent It remains to sum over $i_1$ and $i_2$ as the sums over $l$ and $l'$ are not important (they only
change the constant). We have
\begin{eqnarray*}
\sum_{i_1 \neq i_2} \sum_{l,l'}\left| a_l(i_1) a_l(i_2) a_{l'}(i_1) b_{l'}(i_2) \right| & \leq & 
\sum_{i_1,i_2}  \sum_{l,l'}\left| a_l(i_1) a_l(i_2) a_{l'}(i_1) b_{l'}(i_2) \right| \\
 & \leq & \frac 1 2  \sum_{l,l'} \left( \sum_{i_1, i_2} a_l(i_1)^2 b_{l'}(i_2)^2 + \sum_{i_1, i_2}
a_{l'}(i_1)^2 a_{l}(i_2)^2 \right) \\
 & = & \frac {n^2} 2  \sum_{l,l'}\left( \left< a_l, a_l \right>_n \left< b_{l'}, b_{l'} \right>_n + \left< a_{l'}, a_{l'} \right>_n
\left< a_{l}, a_{l} \right>_n\right)\\
& \leq & \frac{M^2n^2}{K^2}.
\end{eqnarray*}
Clearly, this term behaves like $n^2$. \noindent The other covariances involving three random variables are handled exactly in the
same way.\\

\noindent By combining all the previous bounds, we conclude that
$$A_1 \leq (k_1 + k_2 \ 2^j) \frac{M^2n^2}{K^2},  \quad \mbox{with} \quad 
 k_1 =  224RL^3\|\phi\|^2_\infty, \quad 
k_2 = 32RL^2\|\phi\|^2_\infty.$$

\noindent As a consequence if we write $\bar{C}_{_T}=k_1+k_2$ one gets
$$A_1 \leq \bar{C}_{_T} \frac{M^2}{K^2} \  2^j \ n^2.$$
\hfill
$\Box$ 

\begin{lemma}\label{3ind}
There exists a constant $\tilde{C}_{_T}=\tilde{C}_{_T}(R,L,\|\phi\|_\infty)>0$ such that for any $j \in \mathbb{N}$ $$A_3:=\sum_{i_1 \neq i_2 \neq i_3} \mathbb{C}ov\left( h_j\left( i_1, i_2 \right),
h_j\left( i_1, i_3 \right) \right) \leq \tilde{C}_{_T} \frac{M^2}{K^2} \ \left[ n^3\sum_l\|p_l-q_l\|^2_2+ 2^{\frac{j}{2}} \ n^{\frac{5}{2}}\sum_l\|p_l-q_l\|_2\right].$$
\end{lemma}

\noindent {\it{Proof of Lemma \ref{3ind}}}:\\
\noindent Clearly, the term $A_3$  can be bounded as follows
\begin{eqnarray*}
\lefteqn{
A_3=\sum_{i_1 \neq i_2 \neq i_3} \mathbb{C}ov\left( h_j\left( i_1, i_2 \right),
h_j\left( i_1, i_3 \right) \right)} \\
 & = & \sum_{i_1 \neq i_2 \neq i_3} \sum_{k,k'} \sum_{l,l'}
\mathbb{C}ov\big( \left( a_l(i_1) \phi_{jk}(Y_{i_1}) - b_l(i_1)
\phi_{jk}(Z_{i_1}) \right) \left( a_l(i_2) \phi_{jk}(Y_{i_2}) - b_l(i_2)
\phi_{jk}(Z_{i_2}) \right), \\
 & & \left( a_{l'}(i_1) \phi_{jk'}(Y_{i_1}) - b_{l'}(i_1)
\phi_{jk'}(Z_{i_1}) \right) \left( a_{l'}(i_3) \phi_{jk'}(Y_{i_3}) - b_{l'}
(i_3) \phi_{jk'}(Z_{i_3}) \right) \big) \\
 & = & \sum_{i_1 \neq i_2 \neq i_3} \sum_{k,k'} \sum_{l,l'}
\mathbb{C}ov\left(  a_l(i_1) \phi_{jk}(Y_{i_1}) - b_l(i_1)
\phi_{jk}(Z_{i_1}),
a_{l'}(i_1) \phi_{jk'}(Y_{i_1}) - b_{l'}(i_1)
\phi_{jk'}(Z_{i_1})   \right)\\
&  & \quad \quad \quad \quad \times \mathbb E\left(  a_l(i_2) \phi_{jk}(Y_{i_2}) - b_l(i_2)
\phi_{jk}(Z_{i_2})  \right)
\mathbb E\left( a_{l}(i_3) \phi_{jk'}(Y_{i_3}) - b_{l'}
(i_3) \phi_{jk'}(Z_{i_3})  \right) \\
  & = & \sum_{i_1, i_2, i_3} \sum_{k,k'} \sum_{l,l'}
\mathbb{C}ov\left(  a_l(i_1) \phi_{jk}(Y_{i_1}) - b_l(i_1)
\phi_{jk}(Z_{i_1}),  a_{l'}(i_1) \phi_{jk'}(Y_{i_1}) - b_{l'}(i_1)
\phi_{jk'}(Z_{i_1})  \right) \\
 & & \mathbb E\left( a_l(i_2) \phi_{jk}(Y_{i_2}) - b_l(i_2)
\phi_{jk}(Z_{i_2}) \right)
\mathbb E\left(  a_{l'}(i_3) \phi_{jk'}(Y_{i_3}) - b_{l'}
(i_3) \phi_{jk'}(Z_{i_3}) \right) \\
 & & - \sum_{i_1 = i_2, i_3} \sum_{k,k'} \sum_{l,l'}
\mathbb{C}ov\left( a_l(i_1) \phi_{jk}(Y_{i_1}) - b_l(i_1)
\phi_{jk}(Z_{i_1}) ,
 a_{l'}(i_1) \phi_{jk'}(Y_{i_1}) - b_{l'}(i_1)
\phi_{jk'}(Z_{i_1})  \right) \\
 & & \mathbb E\left( a_l(i_2) \phi_{jk}(Y_{i_2}) - b_l(i_2)
\phi_{jk}(Z_{i_2}) \right)
\mathbb E\left( a_{l'}(i_3) \phi_{jk'}(Y_{i_3}) - b_{l'}
(i_3) \phi_{jk'}(Z_{i_3})  \right) \\
 & & - \sum_{i_1 = i_3, i_2} \sum_{k,k'} \sum_{l,l'}
\mathbb{C}ov\left(  a_l(i_1) \phi_{jk}(Y_{i_1}) - b_l(i_1)
\phi_{jk}(Z_{i_1}) ,
 a_{l'}(i_1) \phi_{jk'}(Y_{i_1}) - b_{l'}(i_1)
\phi_{jk'}(Z_{i_1})   \right) \\
 & & \mathbb E\left( a_l(i_2) \phi_{jk}(Y_{i_2}) - b_l(i_2)
\phi_{jk}(Z_{i_2}) \right)
\mathbb E\left( a_{l'}(i_3) \phi_{jk'}(Y_{i_3}) - b_{l'}
(i_3) \phi_{jk'}(Z_{i_3}) \right) \\
 & & + \sum_{i_1 = i_2 =  i_3} \sum_{k,k'} \sum_{l,l'}
\mathbb{C}ov\left( a_l(i_1) \phi_{jk}(Y_{i_1}) - b_l(i_1)
\phi_{jk}(Z_{i_1}) ,
 a_{l'}(i_1) \phi_{jk'}(Y_{i_1}) - b_{l'}(i_1)
\phi_{jk'}(Z_{i_1}) \right) \\
 & & \mathbb E\left( a_l(i_2) \phi_{jk}(Y_{i_2}) - b_l(i_2)
\phi_{jk}(Z_{i_2}) \right)
\mathbb E\left( a_{l'}(i_3) \phi_{jk'}(Y_{i_3}) - b_{l'}
(i_3) \phi_{jk'}(Z_{i_3}) \right) \\
 & = & A_{31} - A_{32} - A_{33} + A_{34}\\
  & \leq & |A_{31}|  + |A_{32}| + |A_{33}| + |A_{34}|.
\end{eqnarray*}
We will separetely bound each term.

\noindent Let us start with $|A_{31}|$. The first step is to expand the covariance.
\begin{eqnarray*}
|A_{31}| &=& |\sum_{i_1, i_2, i_3} \sum_{k,k'} \sum_{l,l'}
\mathbb{C}ov\left(  a_l(i_1) \phi_{jk}(Y_{i_1}) - b_l(i_1)
\phi_{jk}(Z_{i_1}),  a_{l'}(i_1) \phi_{jk'}(Y_{i_1}) - b_{l'}(i_1)
\phi_{jk'}(Z_{i_1})  \right) \\
& & \mathbb E\left( a_l(i_2) \phi_{jk}(Y_{i_2}) - b_l(i_2)
\phi_{jk}(Z_{i_2}) \right)| \\
& =  & n^2 \ \Big| \sum_{i_1} \sum_{k,k'} \sum_{l,l'}
 \mathbb{C}ov\left( \left( a_l(i_1) \phi_{jk}(Y_{i_1}) - b_l(i_1)
\phi_{jk}(Z_{i_1}) \right),
\left( a_{l'}(i_1) \phi_{jk'}(Y_{i_1}) - b_{l'}(i_1)
\phi_{jk'}(Z_{i_1}) \right)  \right)  \\
 & & (\int \phi_{jk} p_l - \int \phi_{jk} q_l)(\int \phi_{jk'} p_{l'} - \int \phi_{jk'} q_{l'}) \Big| 
\end{eqnarray*}
\begin{eqnarray*}
 & = & n^2 \ \Big| \sum_{i_1} \sum_{k,k'} \sum_{l,l'}
 \big[ \mathbb E\left( a_l(i_1) \phi_{jk}(Y_{i_1}) 
a_{l'}(i_1) \phi_{jk'}(Y_{i_1}) \right) + 
\mathbb E\left( b_l(i_1) \phi_{jk}(Z_{i_1}) 
 b_{l'}(i_1) \phi_{jk'}(Z_{i_1}) \right) \\
 & & 
 - \mathbb E\left( a_l(i_1) \phi_{jk}(Y_{i_1}) \right)
 \mathbb E\left( a_{l'}(i_1) \phi_{jk'}(Y_{i_1}) \right)
 - \mathbb E\left(  b_l(i_1) \phi_{jk}(Z_{i_1})  \right)
 \mathbb E\left( b_{l'}(i_1) \phi_{jk'}(Z_{i_1}) \right) \big]  \\
 & &   (\int \phi_{jk} p_l - \int \phi_{jk} q_l )(\int \phi_{jk'} p_{l'} - \int \phi_{jk'} q_{l'}) \Big|.
\end{eqnarray*}
The first two terms involve only one expectation and can be bounded in the same way. Therefore let us bound the quantity 
$$\left|\sum_{i_1}\sum_{k,k'}\sum_{l,l'}  \mathbb{E}\left(a_l(i_1) \phi_{jk}(Y_{i_1}) 
a_{l'}(i_1)\phi_{jk'}(Y_{i_1}) \right)
(\int \phi_{jk} p_l - \int \phi_{jk} q_l) 
(\int \phi_{jk'} p_{l'} - \int \phi_{jk'} q_{l'}) \right|.$$

\noindent Clearly $\begin{displaystyle}\sum_{i_1}\end{displaystyle}|a_l(i_1)a_{l'}(i_1)|\leq  n \sqrt{\left< a_{l}, a_{l} \right>_n\left< a_{l'}, a_{l'} \right>_n}\leq \frac{M}{K} \ n$.\\
\noindent Since $|\mathbb E\left( \phi_{jk}(Y_{i_1}) \phi_{jk'}(Y_{i_1})
\right)| \leq \begin{displaystyle}\sup_l\end{displaystyle}(\|p_l\|_\infty \vee \|q_l\|_\infty)$, lemma \ref{teclem3} entails that
$$\sum_{k}\sum_{k'  : I_{jk} \cap I_{jk'} \neq \emptyset}
 \left| (\int \phi_{jk} p_l - \int \phi_{jk} q_l)(
\int \phi_{jk'} p_{l'} - \int \phi_{jk'} q_{l'}) \right| \leq 4L^2 \left(\|p_l-q_l\|_2^2+\|p_{l'}-q_{l'}\|_2^2\right).$$

\noindent Then one deduces that for any $1 \leq i_1 \leq n$

\begin{eqnarray*}
&&\sum_{k, k'}\sum_{l, l'}
 |\mathbb E\left( \phi_{jk}(Y_{i_1}) \phi_{jk'}(Y_{i_1})
\right)| \left| (\int \phi_{jk} p_l - \int \phi_{jk} q_l)(
\int \phi_{jk'} p_{l'} - \int \phi_{jk'} q_{l'}) \right| \\
 & \leq & 
 \: 8L^2\sup_l(\|p_l\|_\infty \vee \|q_l\|_\infty) \  \sum_l\left\| p_l -q_l \right\|_2^2.
\end{eqnarray*}

\noindent Hence 
\begin{eqnarray*}
&&\left|\sum_{i_1}\sum_{k,k'}\sum_{l,l'}  \mathbb{E}\left(a_l(i_1) \phi_{jk}(Y_{i_1}) 
a_{l'}(i_1)\phi_{jk'}(Y_{i_1}) \right)
(\int \phi_{jk} p_l - \int \phi_{jk} q_l) 
(\int \phi_{jk'} p_{l'} - \int \phi_{jk'} q_{l'}) \right|\\
&\leq& \frac{8ML^2}{K} \ \sup_l(\|p_l\|_\infty \vee \|q_l\|_\infty) \  \sum_l\left\| p_l -q_l \right\|_2^2 \ n.
\end{eqnarray*}

\noindent Now we come to the last two terms
 which involve two expectations. Let us consider
for example the quantity
\begin{eqnarray*}
\lefteqn{\left|\sum_{k,k'} \mathbb{E}\left( \phi_{jk}(Y_{i_1}) \right)
 \mathbb E\left( \phi_{jk'}(Y_{i_1}) \right)
(\int \phi_{jk} p_l - \int \phi_{jk} q_l) 
(\int \phi_{jk'} p_{l'} - \int \phi_{jk'} q_{l'}) \right|} \\
  & \leq & 
 \sum_{k,k'} \left| \mathbb{E}\left( \phi_{jk}(Y_{i_1}) \right)
 \mathbb E\left( \phi_{jk'}(Y_{i_1}) \right) \right|
\frac 1 2 \left\{ 
\left( \int \phi_{jk} p_l - \int \phi_{jk} q_l  \right)^2
+ \left( \int \phi_{jk'} p_{l'}- \int \phi_{jk'} q_{l'}  \right)^2 \right\}\\
&\leq& \sup_{l,l'}\left[ \sqrt{2L} \left(\|p_l\|_\infty\vee \|q_l\|_\infty\right)^{\frac{1}{2}} \ 2^{-\frac{j}{2}}\sum_k |\mathbb{E}\left( \phi_{jk}(Y_{i_1}) \right)|\sum_{k'} \left( \int \phi_{jk'} p_{l'}- \int \phi_{jk'} q_{l'}  \right)^2\right]\\
&\leq &4\sqrt{2}L^{\frac{5}{2}} \|\phi\|_\infty \left(\sup_l(\|p_l\|_\infty\vee \|q_l\|_\infty)\right)^{\frac{1}{2}}  \sup_{l'}\left\| p_{l'} -q_{l'} 
\right\|_2^2.
\end{eqnarray*}
Last inequalities are obtained by using lemma \ref{teclem3} for any $1 \leq i_1 \leq n.$ Hence 
\begin{eqnarray*}
&&\left|\sum_{i_1}\sum_{k,k'}\sum_{l,l'}  \mathbb{E}\left(a_l(i_1) \phi_{jk}(Y_{i_1})\right) 
\mathbb{E}\left(a_{l'}(i_1)\phi_{jk'}(Y_{i_1}) \right)
(\int \phi_{jk} p_l - \int \phi_{jk} q_l) 
(\int \phi_{jk'} p_{l'} - \int \phi_{jk'} q_{l'}) \right|\\
&\leq& \frac{4\sqrt{2}M}{K}L^{\frac{5}{2}} \|\phi\|_\infty \left(\sup_l(\|p_l\|_\infty\vee \|q_l\|_\infty)\right)^{\frac{1}{2}} \sup_{l'} \left\| p_{l'} -q_{l'} 
\right\|_2^2 \ n.
\end{eqnarray*}

\noindent Therefore the two last bounds entail that
$$\left| A_{31} \right| \leq c_{31} \frac{M^2n^3}{K}  \sum_{l} \left\| p_l -q_l 
\right\|_2^2, \quad \mbox{where} \quad
 c_{31}=4L^2\sqrt{R} \left(2\sqrt{R}+\sqrt{2L}\|\phi\|_\infty\right).$$

\noindent The way to bound $A_{32}$ and $A_{33}$ is trickier.
We have
\begin{eqnarray*}
\left| A_{32} \right| & \leq & \Big| \sum_{l,l'} \sum_{i_1,i_2} \sum_{k,k'}
\left[ a_l(i_1) a_{l'}(i_1) \mathbb{C}ov\left( \phi_{jk}(Y_{i_1}),
\phi_{jk'}(Y_{i_1}) \right) + b_l(i_1) b_{l'}(i_1) \mathbb{C}ov\left( \phi_{jk}(Z_{i_1}), 
\phi_{jk'}(Z_{i_1}) \right) \right]  \\
 & & 
\left[ a_l(i_1) \mathbb E\left( \phi_{jk}(Y_{i_1}) \right) - b_l(i_1) \mathbb E\left(
\phi_{jk}(Z_{i_1}) \right) \right]
\left[ a_{l'}(i_2) \mathbb E\left( \phi_{jk'}(Y_{i_2}) \right) - b_{l'}(Z_{i_2})
 \mathbb E\left(
\phi_{jk'}(Z_{i_2}) \right) \right] \Big| \\
& \leq & \Big| \sum_{l,l'} \sum_{i_1} \sum_{k,k'}
\left[ a_l(i_1) a_{l'}(i_1) \mathbb{C}ov\left( \phi_{jk}(Y_{i_1}),
\phi_{jk'}(Y_{i_1}) \right) + b_l(i_1) b_{l'}(i_1) \mathbb{C}ov\left( \phi_{jk}(Z_{i_1}),
\phi_{jk'}(Z_{i_1}) \right) \right] \\
 & & 
\left[ a_l(i_1) \mathbb E\left( \phi_{jk}(Y_{i_1}) \right) - b_l(i_1) \mathbb E\left(
\phi_{jk}(Z_{i_1}) \right) \right]
\left( n \int \phi_{jk'} p_{l'} - n \int \phi_{jk'} q_{l'} \right) \Big| 
\end{eqnarray*}
\begin{eqnarray*}
 & \leq & \Big| \sum_{l,l'} \sum_{i_1} \sum_{k,k'}
 a_l(i_1) a_{l'}(i_1) \mathbb E\left( \phi_{jk}(Y_{i_1}) \phi_{jk'}(Y_{i_1}) \right)  \\
 & & 
\left( a_l(i_1) \mathbb E\left( \phi_{jk}(Y_{i_1}) \right) - b_l(i_1) \mathbb E\left(
\phi_{jk}(Z_{i_1}) \right) \right)
\left( n \int \phi_{jk'} p_{l'} - n \int \phi_{jk'} q_{l'} \right) \Big| \\
 & &  +  \Big| \sum_{l,l'} \sum_{i_1} \sum_{k,k'}
 a_l(i_1) a_{l'}(i_1) \mathbb E\left( \phi_{jk}(Y_{i_1}) \right)
\mathbb E \left( \phi_{jk'}(Y_{i_1}) \right)  \\
 & & 
\left( a_l(i_1) \mathbb E\left( \phi_{jk}(Y_{i_1}) \right) - b_l(i_1) \mathbb E\left(
\phi_{jk}(Z_{i_1}) \right) \right)
\left( n \int \phi_{jk'} p_{l'} - n \int \phi_{jk'} q_{l'} \right) \Big| \\
 & & + \Big| \sum_{l,l'} \sum_{i_1} \sum_{k,k'}
 b_l(i_1) b_{l'}(i_1) \mathbb E\left( \phi_{jk}(Z_{i_1}) \phi_{jk'}(Z_{i_1}) \right)  \\
 & & 
\left( a_l(i_1) \mathbb E\left( \phi_{jk}(Y_{i_1}) \right) - b_l(i_1) \mathbb E\left(
\phi_{jk}(Z_{i_1}) \right) \right)
\left( n \int \phi_{jk'} p_{l'} - n \int \phi_{jk'} q_{l'} \right) \Big| \\
 & &  +  \Big| \sum_{l,l'} \sum_{i_1} \sum_{k,k'}
 b_l(i_1) b_{l'}(i_1) \mathbb E\left( \phi_{jk}(Z_{i_1}) \right)
\mathbb E \left( \phi_{jk'}(Z_{i_1}) \right)  \\
 & & 
\left( a_l(i_1) \mathbb E\left( \phi_{jk}(Y_{i_1}) \right) - b_l(i_1) \mathbb E\left(
\phi_{jk}(Z_{i_1}) \right) \right)
\left( n \int \phi_{jk'} p_{l'} - n \int \phi_{jk'} q_{l'} \right) \Big|.
\end{eqnarray*}
The calculations are rather lengthy and involve eight terms. But the bright
 side is that the terms can be split into two groups. There are terms
involving two expectations such as
$$
\sum_{l,l'} \sum_{i_1} \sum_{k,k'}
 a_l(i_1) a_{l'}(i_1) \mathbb E\left( \phi_{jk}(Y_{i_1}) \phi_{jk'}(Y_{i_1}) \right) a_l(i_1) \mathbb E\left( \phi_{jk}(Y_{i_1}) \right) 
\left( n \int \phi_{jk'} p_{l'} - n \int \phi_{jk'} q_{l'} \right), $$
and terms involving three expectations such as
$$
\sum_{l,l'} \sum_{i_1} \sum_{k,k'}
 a_l(i_1) a_{l'}(i_1) \mathbb E\left( \phi_{jk}(Y_{i_1}) \right)
\mathbb E \left(
\phi_{jk'}(Y_{i_1}) \right) a_l(i_1) \mathbb E\left( \phi_{jk}(Y_{i_1}) 
\right) 
\left( n \int \phi_{jk'} p_{l'} - n \int \phi_{jk'} q_{l'} \right).
$$
\noindent Still using lemmas \ref{teclem1} and \ref{teclem2}, we have

\begin{eqnarray*}
\lefteqn{\left| \sum_{l,l'} \sum_{i_1} \sum_{k,k'}
 a_l(i_1) a_{l'}(i_1) \mathbb E\left( \phi_{jk}(Y_{i_1}) \phi_{jk'}(Y_{i_1}) \right) a_l(i_1) \mathbb E\left( \phi_{jk}(Y_{i_1}) \right) 
\left( n \int \phi_{jk'} p_{l'} - n \int \phi_{jk'} q_{l'} \right) \right|} \\
& \leq & 8 \sqrt{2}\ L^{\frac 5 2} \|\phi\|_\infty^2
\sup_l \left\| p_l \right\|_\infty^{\frac 1 2}
\sum_{l=1}^M  \left\|p_{l} - q_{l} \right\|_2
\left( \sum_{l=1}^M \left< a_l, a_l \right>_n \right)^{\frac 3 2}  2^{\frac{j}{2}} \ n^{\frac 5 2} \\
& \leq &8 \sqrt{\frac{2M^3}{K^3}}\ L^{\frac 5 2} \|\phi\|_\infty^2
\sup_l \left\| p_l \right\|_\infty^{\frac 1 2}
\sum_{l=1}^M  \left\|p_{l} - q_{l} \right\|_2
 \ 2^{\frac{j}{2}} \ n^{\frac 5 2};\end{eqnarray*}
\begin{eqnarray*}
\lefteqn{\left| \sum_{l,l'} \sum_{i_1} \sum_{k,k'}
 a_l(i_1) a_{l'}(i_1) \mathbb E\left( \phi_{jk}(Y_{i_1}) \phi_{jk'}(Y_{i_1}) \right) b_l(i_1) \mathbb E\left( \phi_{jk}(Z_{i_1}) \right) 
\left( n \int \phi_{jk'} p_{l'} - n \int \phi_{jk'} q_{l'} \right) \right|} \\
& \leq & 4 \sqrt{2}\ L^{\frac 5 2} \|\phi\|_\infty^2
\sup_l \left\| q_l \right\|_\infty^{\frac 1 2}
\sum_{l=1}^M  \left\|p_{l} - q_{l} \right\|_2
\left( \sum_{l=1}^M \left< a_l, a_l \right>_n \right)^{\frac 1 2}  \sum_{l=1}^M \left(\left< a_l, a_l \right>_n + \left< b_l, b_l \right>_n \right) \ 2^{\frac{j}{2}} \ n^{\frac 5 2}  \\
& \leq &8 \sqrt{\frac{2M^3}{K^3}}\ L^{\frac 5 2} \|\phi\|_\infty^2
\sup_l \left\| q_l \right\|_\infty^{\frac 1 2}
\sum_{l=1}^M  \left\|p_{l} - q_{l} \right\|_2
 \ 2^{\frac{j}{2}} \ n^{\frac 5 2};\\
\lefteqn{\left| \sum_{l,l'} \sum_{i_1} \sum_{k,k'}
 b_l(i_1) b_{l'}(i_1) \mathbb E\left( \phi_{jk}(Z_{i_1}) \phi_{jk'}(Z_{i_1}) \right) b_l(i_1) \mathbb E\left( \phi_{jk}(Z_{i_1}) \right) 
\left( n \int \phi_{jk'} p_{l'} - n \int \phi_{jk'} q_{l'} \right) \right|} \\
& \leq & 8 \sqrt{2}\ L^{\frac 5 2} \|\phi\|_\infty^2
\sup_l \left\| q_l \right\|_\infty^{\frac 1 2}
\sum_{l=1}^M  \left\|p_{l} - q_{l} \right\|_2
\left( \sum_{l=1}^M \left< b_l, b_l \right>_n \right)^{\frac 3 2}  2^{\frac{j}{2}} \ n^{\frac 5 2}\\
& \leq &8 \sqrt{\frac{2M^3}{K^3}}\ L^{\frac 5 2} \|\phi\|_\infty^2
\sup_l \left\| q_l \right\|_\infty^{\frac 1 2}
\sum_{l=1}^M  \left\|p_{l} - q_{l} \right\|_2
 \ 2^{\frac{j}{2}} \ n^{\frac 5 2};\\
\lefteqn{\left| \sum_{l,l'} \sum_{i_1} \sum_{k,k'}
 b_l(i_1) b_{l'}(i_1) \mathbb E\left( \phi_{jk}(Z_{i_1}) \phi_{jk'}(Z_{i_1}) \right) a_l(i_1) \mathbb E\left( \phi_{jk}(Y_{i_1}) \right) 
\left( n \int \phi_{jk'} p_{l'} - n \int \phi_{jk'} q_{l'} \right) \right|} \\
& \leq &4 \sqrt{2}\ L^{\frac 5 2} \|\phi\|_\infty^2
\sup_l \left\| p_l \right\|_\infty^{\frac 1 2}
\sum_{l=1}^M  \left\|p_{l} - q_{l} \right\|_2
\left( \sum_{l=1}^M \left< b_l, b_l \right>_n \right)^{\frac 1 2}  \sum_{l=1}^M \left(\left< a_l, a_l \right>_n + \left< b_l, b_l \right>_n \right) \ 2^{\frac{j}{2}} \ n^{\frac 5 2}\\
& \leq &8 \sqrt{\frac{2M^3}{K^3}}\ L^{\frac 5 2} \|\phi\|_\infty^2
\sup_l \left\| p_l \right\|_\infty^{\frac 1 2}
\sum_{l=1}^M  \left\|p_{l} - q_{l} \right\|_2
 \ 2^{\frac{j}{2}} \ n^{\frac 5 2}.
\end{eqnarray*}

\noindent Next we come to the second term. We have
\begin{eqnarray*}
\lefteqn{\left| \sum_{l,l'} \sum_{i_1} \sum_{k,k'}
 a_l(i_1) a_{l'}(i_1) \mathbb E\left( \phi_{jk}(Y_{i_1}) \right)
\mathbb E \left( \phi_{jk'}(Y_{i_1}) \right)  
 a_l(i_1) \mathbb E\left( \phi_{jk}(Y_{i_1}) \right) 
\left( n \int \phi_{jk'} p_{l'} - n \int \phi_{jk'} q_{l'} \right) \right|} & 
& \\
& \leq & 4 \sqrt{2}\|\phi\|_\infty^2 L^{\frac 5 2}
\sup_l \left\| p_l \right\|_\infty^{\frac 1 2}
\sum_{l=1}^M  \left\|p_{l} - q_{l} \right\|_2
\left( \sum_{l=1}^M \left< a_l, a_l \right>_n \right)^{\frac 3 2}  2^{\frac{j}{2}} \ n^{\frac 5 2} \\
& \leq & 4\sqrt{2\frac{M^3}{K^3}}\|\phi\|_\infty^2 L^{\frac 5 2} 
\sup_l \left\| p_l \right\|_\infty^{\frac 1 2}
\sum_{l=1}^M  \left\|p_{l} - q_{l} \right\|_2
\ 2^{\frac{j}{2}} \ n^{\frac 5 2}, 
\end{eqnarray*}
\begin{eqnarray*}
\lefteqn{\left| \sum_{l,l'} \sum_{i_1} \sum_{k,k'}
 a_l(i_1) a_{l'}(i_1) \mathbb E\left( \phi_{jk}(Y_{i_1}) \right)
\mathbb E \left( \phi_{jk'}(Y_{i_1}) \right)  
 b_l(i_1) \mathbb E\left( \phi_{jk}(Z_{i_1}) \right) 
\left( n \int \phi_{jk'} p_{l'} - n \int \phi_{jk'} q_{l'} \right) \right|} & 
& \\
& \leq & 2 \sqrt{2}\|\phi\|_\infty^2 L^{\frac 5 2} 
\sup_l \left\| q_l \right\|_\infty^{\frac 1 2}
\sum_{l=1}^M  \left\|p_{l} - q_{l} \right\|_2
\left( \sum_{l=1}^M \left< a_l, a_l \right>_n \right)^{\frac 1 2}  \sum_{l=1}^M \left(\left< a_l, a_l \right>_n + \left< b_l, b_l \right>_n \right) \ 2^{\frac{j}{2}} \ n^{\frac 5 2}  \\
& \leq & 4\sqrt{2\frac{M^3}{K^3}}\|\phi\|_\infty^2 L^{\frac 5 2} 
\sup_l \left\| q_l \right\|_\infty^{\frac 1 2}
\sum_{l=1}^M  \left\|p_{l} - q_{l} \right\|_2
\ 2^{\frac{j}{2}} \ n^{\frac 5 2}, \\
\lefteqn{\left| \sum_{l,l'} \sum_{i_1} \sum_{k,k'}
 b_l(i_1) b_{l'}(i_1) \mathbb E\left( \phi_{jk}(Z_{i_1}) \right)
\mathbb E \left( \phi_{jk'}(Z_{i_1}) \right)  
 b_l(i_1) \mathbb E\left( \phi_{jk}(Z_{i_1}) \right) 
\left( n \int \phi_{jk'} p_{l'} - n \int \phi_{jk'} q_{l'} \right) \right|} & 
& \\
& \leq & 4 \sqrt{2}\|\phi\|_\infty^2 L^{\frac 5 2} 
\sup_l \left\| q_l \right\|_\infty^{\frac 1 2}
\sum_{l=1}^M  \left\|p_{l} - q_{l} \right\|_2
\left( \sum_{l=1}^M \left< b_l, b_l \right>_n \right)^{\frac 3 2}  2^{\frac{j}{2}} \ n^{\frac 5 2} \\
& \leq & 4\sqrt{2\frac{M^3}{K^3}}\|\phi\|_\infty^2 L^{\frac 5 2} 
\sup_l \left\| q_l \right\|_\infty^{\frac 1 2}
\sum_{l=1}^M  \left\|p_{l} - q_{l} \right\|_2
\ 2^{\frac{j}{2}} \ n^{\frac 5 2}, \\
\lefteqn{\left| \sum_{l,l'} \sum_{i_1} \sum_{k,k'}
 b_l(i_1) b_{l'}(i_1) \mathbb E\left( \phi_{jk}(Z_{i_1}) \right)
\mathbb E \left( \phi_{jk'}(Z_{i_1}) \right)  
 a_l(i_1) \mathbb E\left( \phi_{jk}(Y_{i_1}) \right) 
\left( n \int \phi_{jk'} p_{l'} - n \int \phi_{jk'} q_{l'} \right) \right|} & 
& \\
& \leq & 2 \sqrt{2}\|\phi\|_\infty^2 L^{\frac 5 2} 
\sup_l \left\| p_l \right\|_\infty^{\frac 1 2}
\sum_{l=1}^M  \left\|p_{l} - q_{l} \right\|_2
\left( \sum_{l=1}^M \left< b_l, b_l \right>_n \right)^{\frac 1 2}  \sum_{l=1}^M \left(\left< a_l, a_l \right>_n + \left< b_l, b_l \right>_n \right) \ 2^{\frac{j}{2}} \ n^{\frac 5 2} \\
& \leq & 4  \sqrt{2\frac{M^3}{K^3}}\|\phi\|_\infty^2 L^{\frac 5 2} 
\sup_l \left\| p_l \right\|_\infty^{\frac 1 2}
\sum_{l=1}^M  \left\|p_{l} - q_{l} \right\|_2
 \ 2^{\frac{j}{2}} \ n^{\frac 5 2}. 
\end{eqnarray*}
All these bounds entail that
$$ \left| A_{32} \right| \leq 
c_{32} \left(\frac{M}{K}\right)^{\frac{3}{2}}\ 2^{\frac{j}{2}} \ n^{\frac 5 2} \sum_{l=1}^M \left\| p_l - q_l \right\|_2,$$
\noindent with $c_{32}=48\sqrt{2RL^5}\|\phi\|_\infty^2$. As a consequence, one similarly gets
$$ \left| A_{33} \right| \leq 
c_{33} \left(\frac{M}{K}\right)^{\frac{3}{2}}\ 2^{\frac{j}{2}} \ n^{\frac 5 2} \sum_{l=1}^M \left\| p_l - q_l \right\|_2,$$with $c_{33}=c_{32}$.

\noindent Let us now consider $|A_{34}|$.
\begin{eqnarray*}
|A_{34}|  & \leq  & |\sum_{i_1} \sum_{k,k'} \sum_{l,l'}
\mathbb{C}ov\left( a_l(i_1) \phi_{jk}(Y_{i_1}) - b_l(i_1)
\phi_{jk}(Z_{i_1}) ,
 a_{l'}(i_1) \phi_{jk'}(Y_{i_1}) - b_{l'}(i_1)
\phi_{jk'}(Z_{i_1}) \right) \\
 & & \left. \mathbb E\left( a_l(i_1) \phi_{jk}(Y_{i_1}) - b_l(i_1)
\phi_{jk}(Z_{i_1}) \right)
\mathbb E\left( a_{l'}(i_1) \phi_{jk'}(Y_{i_1}) - b_{l'}
(i_1) \phi_{jk'}(Z_{i_1}) \right)\right|.\\
\end{eqnarray*}

\noindent Once again, we apply the Cauchy-Schwarz inequality,
\begin{eqnarray*}
\left|\sum_{i_1}\sum_{l'}\mathbb E\left( a_{l'}(i_1) \phi_{jk'}(Y_{i_1}) - b_{l'}
(i_1) \phi_{jk'}(Z_{i_1}) \right)\right|&=&n \left|\sum_{l'}\int\phi_{jk'}(p_{l'}-q_{l'})\right|\\
& \leq & n \sum_{l'}\|p_{l'}-q_{l'}\|_2.
\end{eqnarray*}

\noindent According to lemma \ref{teclem2}, we have
for any $1 \leq i_1 \leq n$ and any $l$,
\begin{eqnarray*}
\left|\sum_k\mathbb E\left( a_{l}(i_1) \phi_{jk}(Y_{i_1}) - b_{l}
(i_1) \phi_{jk}(Z_{i_1}) \right)\right|&\leq& \left(|a_l(i_1)| \left|\sum_{k}\int\phi_{jk}f_{i_1}\right|\right)\vee \left(|b_l(i_1)| \left|\sum_{k}\int\phi_{jk}g_{i_1}\right|\right)\\
& \leq & 2L (|a_l(i_1)| \vee |b_l(i_1)|)\ 2^{\frac{j}{2}} \ \|\phi\|_\infty\\
& \leq & 2L   \ \|\phi\|_\infty 2^{\frac{j}{2}} \left(\sum_{l=1}^M \left< a_l, a_l \right>_n\vee \sum_{l=1}^M \left< b_l, b_l \right>_n\right)^{\frac{1}{2}}  \sqrt{n}\\
& \leq & 2L   \sqrt{\frac{M}{K}} \ \|\phi\|_\infty 2^{\frac{j}{2}}\  \sqrt{n}.
\end{eqnarray*}

\noindent According to lemmas \ref{teclem1} and \ref{teclem2}, we have
for any fixed $k$,
\begin{eqnarray*}
&&\sum_{i_1}\sum_{k'}|a_l(i_1)a_{l'}(i_1)|\left(\left| \int\phi_{jk}\phi_{jk'}f_{i_1}\right|+\left| \int\phi_{jk}f_{i_1}\int\phi_{jk'}f_{i_1}\right|\right)\\
& \leq& n \sqrt{\left< a_l, a_l \right>_n\left< a_{l'}, a_{l'} \right>_n} 
\ \left(\sum_{k'}\left| \int\phi_{jk}\phi_{jk'}f_{i_1}\right|+\left| \int\phi_{jk}f_{i_1}\right|\sum_{k'}\left|\int\phi_{jk'}f_{i_1}\right|\right)\\
& \leq &  n \sqrt{\left< a_l, a_l \right>_n\left< a_{l'}, a_{l'} \right>_n}
\left( 4L  \ \sup_l(\|p_l\|_\infty \vee \|q_l\|_\infty) + (2L)^{\frac{3}{2}}
\|\phi\|_\infty\sup_l(\|p_l\|_\infty \vee \|q_l\|_\infty)^{\frac{1}{2}}\right)\\
& \leq &  \frac{M}{K}
\left( 4L  \ \sup_l(\|p_l\|_\infty \vee \|q_l\|_\infty) + (2L)^{\frac{3}{2}}
\|\phi\|_\infty\sup_l(\|p_l\|_\infty \vee \|q_l\|_\infty)^{\frac{1}{2}}\right)\ n
\end{eqnarray*}
and
\begin{eqnarray*}
&&\sum_{i_1}\sum_{k'}|b_l(i_1)b_{l'}(i_1)|\left(\left| \int\phi_{jk}\phi_{jk'}g_{i_1}\right|+\left| \int\phi_{jk}g_{i_1}\int\phi_{jk'}g_{i_1}\right|\right)\\
& \leq& n \sqrt{\left< b_l, b_l \right>_n\left< b_{l'}, b_{l'} \right>_n} 
\ \left(\sum_{k'}\left| \int\phi_{jk}\phi_{jk'}g_{i_1}\right|+\left| \int\phi_{jk}g_{i_1}\right|\sum_{k'}\left|\int\phi_{jk'}g_{i_1}\right|\right)\\
& \leq &  n \sqrt{\left< b_l, b_l \right>_n\left< b_{l'}, b_{l'} \right>_n}
\left( 4L \ \sup_l(\|p_l\|_\infty \vee \|q_l\|_\infty) + (2L)^{\frac{3}{2}}
\|\phi\|_\infty\sup_l(\|p_l\|_\infty \vee \|q_l\|_\infty)^{\frac{1}{2}}\right)\\
& \leq &  \frac{M}{K}
\left( 4L  \ \sup_l(\|p_l\|_\infty \vee \|q_l\|_\infty) + (2L)^{\frac{3}{2}}
\|\phi\|_\infty\sup_l(\|p_l\|_\infty \vee \|q_l\|_\infty)^{\frac{1}{2}}\right)\ n.
\end{eqnarray*}

\noindent Hence,
$$|A_{34}|\leq c_{34} \ 2^{\frac{j}{2}} \left(\frac{M}{K}\right)^{\frac{3}{2}} n^{\frac{5}{2}} \sum_l\|p_l-q_l\|_2,$$
with $c_{34}=4L\|\phi\|_\infty \sqrt{R}\left(4L\sqrt{R}+(2L)^{\frac{3}{2}}\|\phi\|_\infty\right)$.\\

\noindent When we carefully look at the bounds of $A_{3i}$  for $i  \in \{1,2,3,4\}$, we deduce that there exists a $Cste>0$ such that  $$A_3 \leq \tilde{C}_{_T}  \frac{M^2}{K^2}\ \left[ n^3\sum_l\|p_l-q_l\|^2_2+ 2^{\frac{j}{2}}  \ n^{\frac{5}{2}}\sum_l\|p_l-q_l\|_2\right],$$
with $\tilde{C}_{_T}=\begin{displaystyle}\sum_{i=1}^4\end{displaystyle}c_{3i}.$ \hfill
$\Box$ \\

\end{document}